\documentclass[11pt,reqno]{amsart}
\usepackage{fullpage,enumerate}
\usepackage{amsfonts}
\usepackage[all]{xy}
\usepackage{amssymb}
\usepackage{comment}
\usepackage{times}
\usepackage{graphicx}
\usepackage{amsmath,bm}
\usepackage[hidelinks]{hyperref}
\usepackage{tikz-cd}
\usepackage{tikz}
\usepackage{appendix}
\usepackage{color}\usepackage{float}
\usepackage[margin=1.3in,marginparwidth=1.1in,marginparsep=0.1in]{geometry}
\usepackage{marginnote,color}
\usepackage{import}
\usepackage[normalem]{ulem}
\usepackage{graphicx}

\usepackage{mathrsfs}
\usepackage{geometry}
\usepackage{epsfig,amscd,amsthm,mathtools,fourier}
\numberwithin{equation}{section}

\newtheorem{thm}{Theorem}[section]
\newtheorem{cor}[thm]{Corollary}
\newtheorem{lem}[thm]{Lemma}
\newtheorem{prop}[thm]{Proposition}

\newtheorem*{ThmA}{Theorem A}
\newtheorem*{ThmB}{Theorem B}
\newtheorem*{ThmC}{Theorem C}

\theoremstyle{definition}
\newtheorem{ex}[thm]{Example}
\newtheorem{nota}[thm]{Notation}
\newtheorem{defn}[thm]{Definition}

\newtheorem{constr}[thm]{Construction}

\newcommand{\irr}{\mathrm{irr}}

\newcommand{\trop}{\mathrm{trop}}

\newcommand{\w}{\mathrm{w}}

\newcommand{\supp}{{\rm Supp}}
\newcommand{\val}{{\rm val}}

\newcommand{\Spec}{{\rm Spec}}

\newcommand{\FF}{{\mathbb F}}

\newcommand{\RR}{{\mathbb R}}
\newcommand{\ZZ}{{\mathbb Z}}
\newcommand{\PP}{{\mathbb P}}

\newcommand{\GG}{{\mathbb G}}

\newcommand{\QQ}{{\mathbb Q}}

\newcommand{\AAA}{{\mathbb A}}

\newcommand{\CL}{{\mathcal L}}

\newcommand{\CO}{{\mathcal O}}
\newcommand{\CP}{{\mathcal P}}
\newcommand{\CM}{{\mathcal M}}

\newcommand{\CC}{{\mathcal C}}

\newcommand{\CH}{{\mathcal H}}

\newcommand{\ft}{{\mathfrak t}}
\newcommand{\fh}{{\mathfrak h}}

\newcommand{\fw}{{\mathfrak w}}

\newcommand{\tth}{{\rm\mathbf h}}

\usepackage{scalerel}

\newcommand{\subdivision}{{\scaleobj{.5}{\begin{tikzpicture}[x=.52pt,y=.62pt,yscale=1,xscale=1]
\draw   [line width=0.4mm](10,0) -- (-10,0) ;
\draw [line width=0.6mm] (10,0)--(0,20);
\draw (-10,0)--(0,20);
\draw    (-5,10) -- (5,10) ;
\draw    (-5,10) -- (0,0) ;
\draw    (5,10) -- (0,0) ;
\end{tikzpicture}}}}

\newcommand{\Star}{{\rm Star}}

\def\:{\colon} 
\def\val{\nu}

\newcommand{\cC}{{\mathcal{C}}}

\newcommand{\m}{\mathfrak w}

\theoremstyle{remark}
\newtheorem{rem}[thm]{Remark}

\makeatletter
\DeclareRobustCommand{\cev}[1]{
  {\mathpalette\do@cev{#1}}
}
\newcommand{\do@cev}[2]{
  \vbox{\offinterlineskip
    \sbox\z@{$\m@th#1 x$}
    \ialign{##\cr
      \hidewidth\reflectbox{$\m@th#1\vec{}\mkern4mu$}\hidewidth\cr
      \noalign{\kern-\ht\z@}
      $\m@th#1#2$\cr
    }
  }
}
\makeatother

\title{The irreducibility of Hurwitz spaces and Severi varieties on toric surfaces}
\author{Karl Christ, Xiang He, and Ilya Tyomkin}
\thanks{KC was partially supported by NSF FRG grant DMS–2053261, a Minerva Short-Term Research Grant, and the Center for Advanced Studies at BGU. XH is supported by the NSFC grant 12301057, the National Key R\&D Program of China 2022YFA1007100, and the ERC Consolidator Grant 770922 - BirNonArchGeom. }

\address[Christ]{Dipartimento di Matematica\\
    Università di Torino\\Via Carlo Alberto 10 \\10123 Turin\\  Italy }\email{karl.christ@unito.it}

\address[He]{Yau Mathematical Sciences Center\\ Shuangqing Complex Building, Tsinghua University\\ Haidian District, Beijing\\ 100084\\ China}\email{xianghe@mail.tsinghua.edu.cn}

\address[Tyomkin]{Department of Mathematics\\
    Ben-Gurion University of the Negev\\P.O.Box 653 \\Be'er Sheva\\ 84105\\  Israel}\email{tyomkin@math.bgu.ac.il}

\begin{document}
\tikzset{every picture/.style={line width=0.75pt}} 

\begin{abstract}
    In 1969, Fulton introduced classical Hurwitz spaces parametrizing simple $d$-sheeted coverings of the projective line in the algebro-geometric setting. He established the irreducibility of these spaces under the assumption that the characteristic of the ground field is greater than $d$, but the irreducibility problem in smaller characteristics remained open. We resolve this problem in the current paper and prove that the classical Hurwitz spaces are irreducible over any algebraically closed field. Along the way, we establish the irreducibility of Severi varieties in arbitrary characteristic for a rich class of toric surfaces, including all classical toric surfaces. Our approach to the irreducibility problems comes from tropical geometry, and the paper contains two more results of independent interest---a lifting result for parametrized tropical curves and a strong connectedness property of the moduli spaces of parametrized tropical curves.  
\end{abstract}

\maketitle
    
\setcounter{tocdepth}{1}
\tableofcontents

\section{Introduction}

In this paper, we establish the irreducibility of the classical Hurwitz spaces over algebraically closed fields in {\em arbitrary characteristic}. The classical Hurwitz spaces $H_{g,d}$ classify simple $d$-sheeted coverings of the projective line by smooth genus-$g$ curves. Here, {\em simple} means that over each branch point, there exists a unique ramification point, and the ramification at this point is as simple as possible, i.e., if the characteristic is different from 2, then the ramification index is 2, and if the characteristic is 2, then the local degree at the ramification point is 2 and the sheaf of relative differentials has length 2.

Given the central place Hurwitz spaces occupy in the theory of algebraic curves, establishing a property as basic as irreducibility should have far-reaching consequences. For example, it allows one to remove all assumptions on the characteristic in the recent advances in the Brill-Noether theory of general $k$-gonal curves \cite{LLV25, CJLLV}, as well as in related irreducibility results for the space of certain linearly normal curves in $\PP^r$ \cite{Keem}; see \cite[p.4]{Keem}, \cite[Remark 1]{LLV25}, and \cite[Remark 1.5]{CJLLV}.

Classical Hurwitz spaces were first introduced in the complex-analytic setting by Hurwitz in 1891 \cite{Hur91}. Hurwitz showed that these spaces admit a natural structure of a complex manifold. Then, using a monodromy argument going back to Clebsch and Lüroth \cite{Lur71, Cle73}, he proved that the spaces $H_{g,d}$ are connected.

In 1921, Severi proposed several proofs of the irreducibility of the moduli spaces $M_g$ of compact Riemann surfaces of genus $g$. One of the proofs was based on the irreducibility of Hurwitz spaces, using the observation that any compact Riemann surface $C$ of genus $g$ admits a meromorphic function $f$ of degree $d$ such that the covering $(C,f)$ is simple for any $d>g+1$ \cite{Severi}. This circle of ideas goes back to Klein, who already in 1882 noticed that the connectedness of $H_{g,d}$ for large $d$ implies the connectedness of the moduli spaces $M_g$, even before these spaces were formally defined \cite{Klein}.

An algebraic definition of classical Hurwitz spaces is due to Fulton in 1969 \cite{Ful69}, who proved the representability of the corresponding functor\footnote{Fulton works over an arbitrary base scheme, and his notion of simplicity in characteristic 2 is stronger than the one mentioned above. As a result, there are no simple coverings in Fulton's sense in characteristic 2.} for $d>2$ by a Noetherian scheme \cite[Corollary~6.4]{Ful69}. Fulton showed that the spaces $H_{g,d}$ are smooth and irreducible if the characteristic of the ground field is at least $d+1$, and deduced from this the irreducibility of the moduli spaces $\CM_g$ of smooth projective curves if the characteristic is greater than $g+1$. Independently, and using different methods, Deligne and Mumford proved the irreducibility of the moduli spaces $\CM_g$ in arbitrary characteristic in 1969 \cite{DM69}.

In this paper, we prove the irreducibility of classical Hurwitz spaces $H_{g,d}$ in arbitrary characteristic. To do so, we consider the Severi varieties of $\PP^1\times\PP^1$ that parametrize irreducible integral genus-$g$ curves of bidegree $(c,d)$ for $c\gg 1$. We prove the irreducibility of the Severi varieties using tools coming from tropical geometry and show that they dominate the spaces $H_{g,d}$, which implies the desired irreducibility result for the classical Hurwitz spaces.

The methods we develop allow us to establish the irreducibility of Severi varieties in arbitrary characteristic for a much richer class of polarized toric surfaces (including all classical toric surfaces such as $\PP^2$, Hirzebruch surfaces, and toric del Pezzo surfaces), providing a far-reaching generalization of our irreducibility result for planar curves \cite{CHT23}. Our approach is different from the one in \cite{CHT23} and relies heavily on our recent theory of tropicalizations of families of parametrized curves over arbitrary bases developed in \cite{CHT24a}.

Finally, let us mention that for a long time, $\PP^2$ was the only surface for which the Severi varieties were known to be irreducible, a result due to Harris \cite{Har86}. Recently, however, the question has received considerable attention for other surfaces. The irreducibility was proved by the third author for Hirzebruch surfaces in characteristic zero in \cite{Tyo07}. Further (ir)reducibility results with restrictions on the genus and the polarization have been obtained for toric surfaces, K3 surfaces, abelian surfaces, etc.; see, e.g., \cite{CC99, Tes09, Tyo14, Bou16, CFGK17, Zar22, LT23, BLC23, CGY23}.

\subsection{Main Results}

Our first main result is the following theorem.

\begin{ThmA}
    The classical Hurwitz spaces $H_{g,d}$ over any algebraically closed field are non-empty and irreducible for all $g\ge 0$ and $d>1$.
\end{ThmA}

To achieve the irreducibility of the classical Hurwitz spaces, we first prove the irreducibility of Severi varieties for a rich class of toric surfaces. This is our second main result.

\begin{ThmB}
    Let $\Delta\subset\RR^2$ be an admissible lattice polygon, and $(S_\Delta,\CL_\Delta)$ the corresponding polarized toric surface. Then, for any $0 \leq g \leq |\Delta^\circ \cap \ZZ^2|$, the Severi variety $V_{g, \Delta}^{\irr}$ parametrizing integral curves of genus $g$ in the linear system $|\CL_\Delta|$ is non-empty and irreducible over any algebraically closed field.
\end{ThmB}

Admissible lattice polygons form a subclass of the so-called $h$-transverse polygons. Their definition (Definition~\ref{def:admissible}) is rather technical; thus, we do not include it in the introduction. For a slightly more restrictive class of {\em very admissible} polygons (Definition~\ref{def:veryadmis}), which still includes the polygons of all classical polarized toric surfaces, we can prove more. We show in Theorem~\ref{thm:zariski} not only that the Severi varieties associated to very admissible polygons are irreducible, but that their general points correspond to at-worst-nodal curves. The latter is a generalization of a remarkable theorem of Zariski \cite{Zar82}; see also \cite{AC81} for a closely related result. As an immediate corollary (Corollary~\ref{cor:adjacencies}), we obtain the following sequence of adjacencies: $\overline V^\irr_{0,\Delta}\subset \overline V^\irr_{1,\Delta}\subset\dots\subset \overline V^\irr_{g,\Delta}$ for very admissible polygons $\Delta$.

Our third main result is purely tropical. Together with the theory of tropicalizations of families of parametrized curves developed in \cite{CHT24a}, and with the lifting result proved in the current paper (Theorem~\ref{thm:existence and uniqueness}), it implies the irreducibility of Severi varieties. Roughly speaking, this result says that on the tropical side, the moduli spaces of parametrized tropical curves of certain degrees satisfy a strong connectedness property. Namely, any two maximal strata classifying tropical curves satisfying vertically stretched point constraints are connected through a sequence of codimension-one strata classifying weightless parametrized tropical curves, all of whose vertices are 3-valent except for a unique 4-valent vertex. Curves satisfying vertically stretched point constraints are called {\em simple stretched floor decomposed (ssfd) curves} (\S~\ref{subsubsec:floor decomposed}). The maximal strata classifying such curves are called {\em ssfd strata} (\S~\ref{subsubsec:moduli of parametrized tropical curves}) and the codimension-one strata as above are called {\em simple walls} (Definition~\ref{def:nice}). Finally, if two strata are connected through a sequence of simple walls, then we say that these strata are {\em sw-equivalent} (Definition~\ref{def:swequiv}).

\begin{ThmC}
    Let $\Delta\subset \RR^2$ be an $h$-transverse polygon and $\nabla$ a tropical degree dual to $\Delta$. Suppose the non-vertical slopes in $\nabla$ are all primitive and $\nabla$ contains either $(0,1)$ or $(0,-1)$. Then any two ssfd strata in $M_{g,\nabla}^\trop$ are sw-equivalent.
\end{ThmC}

\subsection{The ideas of the proofs}

To prove Theorem~A, we pick $c\gg 1$ relatively prime to $d$ and consider the Severi variety parametrizing integral curves in $\PP^1\times\PP^1$ of geometric genus $g$ and bidegree $(c,d)$. Since $c\gg 1$, any $d$-sheeted covering $C\to\PP^1$ of genus $g$ admits a morphism $C\to \PP^1$ of degree $c$, and since $c$ and $d$ are relatively prime, the induced map to $\PP^1\times\PP^1$ is birational onto its image, and therefore defines a point in the Severi variety. Furthermore, it is not difficult to check that if $[D]$ is a general point of the Severi variety, $C$ is the normalization of $D$, and $f$ is the composition of the normalization map with the projection to the first factor, then $(C,f)$ is a simple covering. This shows that the space $H_{g,d}$ is dominated by an alteration of the Severi variety over which the universal curve is equinormalizable. If the Severi variety is irreducible, the alteration can be chosen to be irreducible too, implying the irreducibility of $H_{g,d}$. In other words, this reduces Theorem~A to Theorem~B.

Next, let us explain our approach to Theorem~B. Since Severi varieties are defined over the algebraic closures of primary fields, and since the geometric irreducibility property is stable under field extensions, we may assume that the base field is the algebraic closure $K$ of a complete discretely valued field, which allows us to employ tropical techniques. To an irreducible component $V$ of the Severi variety $V_{g,\Delta}^\irr$, we associate the set $\Sigma_V$ of tropicalizations of curves parametrized by the $K$-points of $V$ that intersect the boundary divisor of $\PP^1\times\PP^1$ transversely. It follows from \cite{CHT24a} that $\Sigma_V$ is the set of rational points of a polyhedral subcomplex of the moduli space $M_{g,\nabla}^\trop$ of parametrized tropical curves of genus $g$ and tropical degree $\nabla$ dual to the polygon $\Delta$. The proof of Theorem~B contains three main ingredients.

The first ingredient is the following balancing property of the polyhedral complex $\overline{\Sigma}_V$ that has been established in \cite{CHT24a}: (a) if $M_{[\Theta]}\subset M_{g,\nabla}^\trop$ is a nice stratum (cf. Definition~\ref{def:nice}) and $\dim(\overline{\Sigma}_V \cap M_{[\Theta]})=\dim M_{[\Theta]}$, then $M_{[\Theta]}\subseteq \overline{\Sigma}_V$, and (b) if $M_{[\Theta]}$ is a simple wall and $\overline{\Sigma}_V$ contains one of the nice strata adjacent to it, then $\overline{\Sigma}_V$ contains the other two as well (Lemma~\ref{lem:surjectivity into nice cones}). Note that ssfd strata are particular cases of nice strata. Therefore, the results of {\em loc. cit.} apply to them.

The second ingredient is Theorem~C, asserting that the moduli space $M_{g,\nabla}^\trop$ satisfies the following connectedness property under suitable assumptions on the tropical degree $\nabla$: any two ssfd strata are connected through a sequence of simple walls. Combined with the first ingredient, this shows that if $\overline{\Sigma}_V$ contains one ssfd stratum, then it contains all such strata. It is not difficult to prove that for any $V$, the complex $\overline{\Sigma}_V$ contains at least one (and hence all) such strata. The latter follows easily from the results of Brugall\'e and Mikhalkin \cite{BM08} by imposing $\dim(V)$ point constraints with vertically stretched tropicalizations on the curves parametrized by $V$.

The third ingredient is the lifting result (Theorem~\ref{thm:existence and uniqueness}) that asserts that if a nice stratum contains a simple parametrized tropical curve of multiplicity one, then it is contained in $\overline{\Sigma}_{V_0}$ for a {\em unique} irreducible component $V_0$ of the Severi variety $V_{g,\Delta}^\irr$. The class of admissible polygons $\Delta$ is precisely the class for which the reduced dual degree $\nabla$ satisfies the assumptions needed for Theorem~C to apply and for which there exists a multiplicity-one ssfd curve of genus $g$. Therefore, for polarized toric surfaces associated to such polygons, we manage to prove the irreducibility of Severi varieties by combining the three ingredients described above.

It is not surprising that multiplicity-one strata satisfy the unique lifting property. Indeed, in his celebrated paper \cite{Mik05}, Mikhalkin proved that the degree of the Severi variety $V_{g,\Delta}^\irr$ is equal to the number of parametrized tropical curves of genus $g$ and degree $\nabla$ satisfying a given tropically general point constraint, and counted with combinatorially defined (Mikhalkin's) multiplicities. In fact, Mikhalkin proved more. His Correspondence Theorem provides a finite-to-one correspondence between algebraic curves in $V_{g,\Delta}^\irr$ satisfying general point constraints and parametrized tropical curves in $M_{g,\nabla}^\trop$ satisfying the tropicalizations of these constraints, as long as the tropicalizations are in tropically general position. Mikhalkin's multiplicity is nothing but the number of algebraic lifts (preimages) of the tropical curve under Mikhalkin's correspondence. In particular, multiplicity-one curves satisfying a given tropically general point constraint admit unique lifts, which explains the intuition behind Theorem~\ref{thm:existence and uniqueness}.

Finally, let us say a few words about the proof of Theorem~C. The advantage of working with ssfd tropical curves is that their deformations are very easy to describe; see \S~\ref{subsec:basic moves}. Our proof proceeds in several steps. First, we show that any ssfd stratum is sw-equivalent to a stratum containing an ssfd curve without self-intersections (Lemma~\ref{lem:self intersection}). Second, to any such stratum, we associate a {\em multiplicity sequence} and show that there exists a unique minimal multiplicity sequence, and that any sw-equivalence class contains a stratum whose multiplicity sequence is the minimal one (Lemma~\ref{lem:minmultseq}). Finally, we show that any two strata whose multiplicity sequences are minimal are sw-equivalent. All steps use two very explicit deformations (Constructions~\ref{constr:moving right} and \ref{constr:moving up}) of ssfd curves that produce the desired strata in a given sw-equivalence class.

\subsection{Reading guide}

In Section~\ref{sec:prelim}, we set up the notation and summarize the specific notions and constructions from algebraic and tropical geometries used throughout the paper. While we provide necessary references to the literature for formal definitions, this section assumes some background knowledge from tropical geometry and is intended primarily as a brisk reminder for the convenience of the reader.

In Section~\ref{sec:lifting}, we prove our lifting result (Theorem~\ref{thm:existence and uniqueness}) that plays an important role in the proof of irreducibility of Severi varieties. The technical core of this section is Lemma~\ref{lem:lifting}.

In Section~\ref{sec:algorithm}, we prove the strong connectedness property of the tropical moduli spaces of parametrized tropical curves of certain degrees (Theorem~C) -- another important ingredient in our approach to the irreducibility problems. The proof of this result is of a combinatorial nature, and the reader interested only in the algebro-geometric irreducibility results can skip it and use Theorem~C as a black box. 

Section~\ref{sec:Severi} is devoted to the proof of irreducibility of Severi varieties on polarized toric surfaces associated to admissible lattice polygons (Theorem~B). We begin this section by defining (very) admissible polygons, and then prove Theorem~\ref{thm:stongB}, which is slightly stronger than Theorem~B. After that we provide a sufficient condition for being very admissible (Proposition~\ref{prop:example of admissible}), which is satisfied by all classical toric surfaces (Example~\ref{ex:very admissible polygons}), and establish a generalization of Zariski's Theorem (Theorem~\ref{thm:zariski}). Again, the reader interested only in the irreducibility results can safely skip these parts and proceed to the following section.

Finally, in Section~\ref{sec:hurwitz}, we define the notion of simple coverings over algebraically closed fields, which differs from Fulton's definition in characteristic 2, and prove the irreducibility of the Hurwitz stack $\CH_{g,d}$ of all $d$-sheeted genus-$g$ coverings (Theorem~\ref{thm:irreducibility of degree d cover}) and of the classical Hurwitz space $H_{g,d}$ of simple $d$-sheeted genus-$g$ coverings (Theorem~A) over algebraically closed fields in all characteristics.

\subsection*{Acknowledgments} 

We would like to thank Eric Larson, Hannah Larson, and Isabel Vogt for helpful discussions and for their interest in our work on the irreducibility of Severi varieties. In particular, they drew our attention to the fact that the irreducibility problem for Hurwitz spaces was open in small characteristic. We also thank the anonymous reviewer for their valuable comments and suggestions.

\section{Preliminaries}\label{sec:prelim}

\subsection{Notation and conventions}\label{sec:not} 

Unless otherwise stated, all tropical curves are connected, and all lattice polygons are convex. 

The graphs we consider have half-edges, called {\em legs}. For a graph $\GG$, we denote the set of vertices by $V(\GG)$, of edges by $E(\GG)$, and of legs by $L(\GG)$. We set $\overline{E}(\GG):=E(\GG)\cup L(\GG)$. For $e\in \overline E(\GG)$, we write $\vec e$ to indicate a choice of an orientation on $e$. By convention, the legs will always be oriented away from the adjacent vertex, and the edges will be considered with both possible orientations. For $\vec{e}\in\overline{E}(\GG)$, we write $\ft(\vec e)$ and $\fh(\vec e)$ to denote the tail and the head of $\vec e$, respectively. For $v\in V(\GG)$, the {\em star} of $v$, i.e., the set of oriented edges and legs having $v$ as their tail, is denoted by $\Star(v)$. In particular, $\Star(v)$ contains a pair of oriented edges associated to any loop adjacent to $v$. We call the size of the set $\Star(v)$ the {\em valence} of $v$.

Throughout this paper, $K$ denotes the algebraic closure of a complete discretely valued field $F$ with an algebraically closed residue field $\widetilde F$. We denote the ring of integers of $K$ by $K^0$, its maximal ideal by $K^{00}$, and the residue field by $\widetilde K=\widetilde F$. The valuation is denoted by $\nu\colon K\rightarrow\mathbb R\cup\{\infty\}$. 

For a toric surface $S$, the lattice of characters is denoted by $M$ and the lattice of cocharacters by $N$. The boundary divisor of $S$ is denoted by $\partial S$, the dense open orbit by $T$, and the monomial functions by $x^m$, for $m\in M$. For a lattice polygon $\Delta\subset M_\RR:=M\otimes\RR$, we denote by $(S_\Delta, \CL_\Delta)$ the polarized toric surface associated to $\Delta$, and by $\Delta^\circ$ and $\partial\Delta$ the interior and the boundary of $\Delta$, respectively.

\subsection{Severi varieties}

For a polarized projective surface $(S,\CL)$ over an algebraically closed field $k$, the {\em Severi variety} $V_{g,\CL}^{\irr}\subset |\CL|$ is the locus of integral curves of geometric genus $g$ that contain no singular points of $S$. By \cite[Lemma~2.6]{CHT23}, $V_{g,\CL}^{\irr}$ is a locally closed subvariety of $|\CL|$. We denote its closure by $\overline V^\irr_{g,\CL}$. If $(S,\CL)=(S_\Delta, \CL_\Delta)$, then we set $V^\irr_{g,\Delta}:=V^\irr_{g,\CL_\Delta}$. By \cite[Proposition~2.7]{CHT23}, $V^\irr_{g,\Delta}$ is either empty or equidimensional of dimension \[\dim(V^\irr_{g,\Delta}) = -\CL_\Delta\cdot K_{S_\Delta}+g-1,\] and any general point $[C]\in V^\irr_{g,\Delta}$ represents a curve $C$ that intersects the boundary divisor transversely; see, e.g., \cite[Proposition~2.7]{CHT23}. Notice, however, that $C$ is not necessarily nodal; see, e.g., \cite[\S~4.1]{Tyo13}.

\subsection{Algebraic curves}\label{subsec: parametrized curves}

A {\em family of curves with marked points} is a flat, proper morphism $\cC\to B$ of finite presentation and relative dimension one, equipped with a tuple of disjoint sections $\sigma_\bullet\:B\to \cC$, contained in the smooth locus of the family. Such a family is called {\em prestable} (resp.\ {\em (semi)stable}) if its geometric fibers have at-worst-nodal singularities (resp.\ are (semi)stable). It is called {\em split} if all its fibers are split, i.e., for any $b\in B$, the irreducible components of the normalization of $\cC_b$ are smooth, geometrically irreducible, and the preimages of the nodes in the normalization are defined over the residue field $k(b)$ (cf.\ \cite[\S~2.21,2.22]{dJ96}). Given an open immersion $B\subset Z$, a {\em model} of $(\cC, \sigma_\bullet)$ over $Z$ is a family of curves with marked points over $Z$ whose restriction to $B$ is the given family $(\cC, \sigma_\bullet)$.

A {\em parametrized curve} in a toric surface $S$ over a field $k$ is a tuple $\left(C,\sigma_\bullet, f\: C \to S\right)$, where $(C,\sigma_\bullet)$ is a smooth projective curve with marked points and $f\: C \to S$ is a morphism such that no zero-dimensional orbit of $S$ is contained in $f(C)$ and the pullback $f^{-1}(\partial S)$ of the boundary divisor is contained in the union of marked points $\bigcup_i\sigma_i$. A {\em family of parametrized curves} $f\colon \cC \rightarrow S$ over $k$ is a family of smooth marked curves $(\cC \to B, \sigma_{\bullet})$ equipped with a morphism $f\colon \cC \rightarrow S$, such that for any geometric point $b\in B$, the restriction $\cC_b \to S$ is a parametrized curve.

\subsection{Tropical curves}\label{subsubsec:abstract and parametrized trop curve}

We follow the notation of \cite{CHT22,CHT23,CHT24a}. Recall that a {\em tropical curve} $\Gamma$ is a finite graph $\GG$ with ordered legs equipped with a pair of functions: the {\em length function} $\ell\colon E(\GG)\rightarrow \RR_{>0}$ and the {\em weight (or genus) function} $g\colon V(\GG)\rightarrow \ZZ_{\geq 0}$. We extend the length function to $\overline{E}(\GG)$ by setting $\ell(l):=\infty$ for any leg $l\in L(\GG)$. If the weight function is identically zero, then we say that the tropical curve is {\em weightless}. The {\em genus} of a tropical curve $\Gamma$ is defined to be $g(\Gamma) = g(\GG) :=1-\chi(\GG)+\sum_{v\in V(\GG)}g(v)$, where $\chi(\GG):=b_0(\GG)-b_1(\GG)$ is the Euler characteristic of $\GG$. A tropical curve $\Gamma$ is said to be {\em stable} if for any vertex, twice its weight plus its valence is at least 3. It is often convenient to consider tropical curves as polyhedral complexes by identifying their edges and legs with closed (semi-)bounded intervals of the corresponding lengths. In this geometric realization, we denote the interior of the interval $e\in \overline{E}(\GG)$ by $e^\circ$. 

Let $N$ be a lattice. Recall that a {\em parametrized tropical curve} is a balanced piecewise integral affine map $h\colon \Gamma\rightarrow N_\mathbb R$ from a tropical curve $\Gamma$ (considered as a polyhedral complex) to the vector space $N_\RR=N\otimes_\ZZ\RR$, where {\em balanced} means that $\sum_{\vec e\in\Star(v)}\frac{\partial h}{\partial \vec e}=0$ for any vertex $v\in V(\Gamma)$. A parametrized tropical curve $(\Gamma,h)$ is {\em stable} if $\Gamma$ is stable. The integral length of the slope $\frac{\partial h}{\partial \vec e}\in N$ is called the {\em multiplicity} of $h$ along $\vec e$, and is denoted by $\m(h,e)$ (or simply $\m(e)$). The {\em degree} $\nabla$ of $(\Gamma,h)$ is the sequence of non-zero slopes $\left(\frac{\partial h}{\partial \vec l_i}\right)$ along the legs of $\Gamma$. We say that a degree $\nabla$ is {\em dual} to a lattice polygon $\Delta\subset M_\RR$ if each slope $\frac{\partial h}{\partial \vec l_j}$ is a multiple of an outer normal of $\Delta$ and the total multiplicity of slopes that correspond to a given side is equal to the integral length of this side. A degree is called {\em reduced} if all vectors in $\nabla$ are primitive.

A parametrized tropical curve $(\Gamma,h)$ is called \textit{simple} (\cite[Definition 4.2]{Mik05}) if it satisfies the following conditions:
\begin{enumerate}
    \item $\Gamma$ is weightless, 3-valent, and $h$ is an immersion;
    \item for any $q\in N_\RR$ the inverse image $h^{-1}(q)$ consists of at most two points;
    \item if $a,b\in \Gamma$ are such that $a\neq b$ and $h(a)=h(b)$, then neither $a$ nor $b$ can be a vertex of $\Gamma$.
\end{enumerate}
A point $q\in N_\RR$ contained in the image $h(\Gamma)$ of a simple parametrized curve $(\Gamma, h)$ is called a \textit{self-intersection point} if $h^{-1}(q)$ consists of two points. 

Next, we recall the notion of multiplicity of $(\Gamma,h)$. Notice that if $v\in V(\Gamma)$ is a 3-valent vertex and $\vec{e_1},\vec{e_2}\in \Star(v)$ a pair of distinct edges, then, by the balancing condition, $\left|\det \left(\frac{\partial h}{\partial \vec e_1},\frac{\partial h}{\partial \vec e_2} \right)\right|$ depends only on the vertex $v$ and not on the choice of the pair of edges. This quantity is called the {\em multiplicity} of the curve $(\Gamma,h)$ at the vertex $v$, and is denoted by ${\rm mult}(v)$ (cf.\ \cite[Definition~2.16]{Mik05}). If $(\Gamma, h)$ is simple, then {\em Mikhalkin's multiplicity} (or simply {\em the multiplicity}) of $(\Gamma, h)$ is defined to be the product of its vertex multiplicities $\prod_{v\in V(\Gamma)}{\rm mult}(v)$ (cf.\ \cite[Definition~4.15]{Mik05}).

\begin{defn} \label{def:unimodular}
     A curve $(\Gamma,h)$ is called \textit{unimodular} if it is simple and 
     \[\left|\det \left(\frac{\partial h}{\partial \vec e_1},\frac{\partial h}{\partial \vec e_2}\right)\right|=1\] 
     for any pair $e_1,e_2\in \overline E(\Gamma)$ such that $e_1\ne e_2$ and $h(e_1)\cap h(e_2)\neq \emptyset$. In particular, any unimodular tropical curve has multiplicity one.
\end{defn}

Finally, recall that to a parametrized tropical curve of degree $\nabla$ dual to $\Delta$, one can associate the so-called {\em dual subdivision} $\subdivision$ of the lattice polygon $\Delta$; see \cite[Proposition 3.11]{Mik05} and Figure~\ref{fig:comtype} for an illustration. If $(\Gamma,h)$ is simple, then the dual subdivision consists of triangles and parallelograms, where the triangles correspond to the vertices of $\Gamma$ and the parallelograms to the self-intersection points. The multiplicities of the edges and legs correspond to the integral lengths of the corresponding edges of the dual subdivision, and the vertex multiplicity of a 3-valent vertex $v$ to twice the Euclidean area of the triangle dual to $v$. In particular, in the dual language, the multiplicity-one condition means that all triangles in the dual subdivision are primitive, i.e., contain no integral points except the vertices, and unimodularity means that both the triangles and the parallelograms are primitive; see Figure~\ref{fig:comtype}. Note that by Pick's formula, a convex lattice $n$-gon in $M_\RR$ is primitive if and only if its Euclidean area is the minimal possible, i.e., $\frac{n-2}{2}$.

\subsection{Floor decomposed tropical curves and $h$-transverse polygons.}\label{subsubsec:floor decomposed}

A particularly convenient combinatorial type of tropical curves was introduced by Brugall\'e and Mikhalkin in \cite{BM08} in the case $N=\ZZ^2$. A {\em floor decomposed tropical curve} is a parametrized tropical curve $h \colon \Gamma \to \RR^2$, in which the $x$-coordinate of the slope of each edge/leg is either 0 or $\pm 1$.

\begin{figure}[ht]
    \begin{tikzpicture}[x=0.8pt,y=0.8pt,yscale=-.6,xscale=.6]
        \import{./}{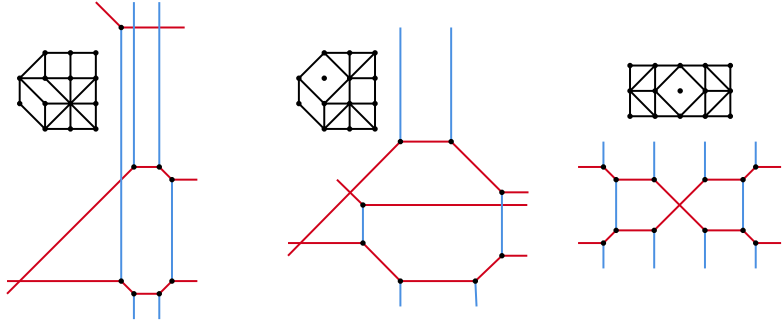}
    \end{tikzpicture}   
    \caption{Three simple floor decomposed curves of multiplicity one and their dual subdivisions. The floors are in red and the elevators in blue. Only the left curve is ssfd and unimodular. The left and the middle curves have different numbers of self-intersection points despite having the same combinatorial type. The right curve does not belong to an ssfd stratum.}
    \label{fig:comtype}
\end{figure}

A non-contracted edge/leg whose slope has $x$-coordinate 0 is called an {\em elevator}. If a connected component of the graph obtained from $\Gamma$ by removing the interiors of all elevators contains a non-contracted edge/leg, it is called a {\em floor}. See \cite[\S~3.2]{CHT23} or \cite{BM08} for details and Figure~\ref{fig:comtype} for an example. 

The existence of a floor decomposed curve depends on the specified degree. Recall that a lattice polygon $\Delta\subset \RR^2$ is called $h${\em -transverse} if its sides have integral or infinite slopes. Equivalently, $\Delta$ is $h$-transverse if and only if for any $r \in \mathbb Z$, the line $y = r$ is either disjoint from $\Delta$ or contains a side of $\Delta$ or intersects its boundary only in lattice points. By \cite{BM08}, $\Delta$ is $h$-transverse if and only if there is a floor decomposed tropical curve of the reduced degree dual to $\Delta$. 

\begin{defn}\label{def:stretched}
    Let $(\Gamma,h)$ be a floor decomposed curve. For $v\in V(\Gamma)$, set $(x_v,y_v):=h(v)\in\RR^2$. 
    \begin{enumerate}
        \item The {\em width} of $(\Gamma,h)$ is defined to be $\w(\Gamma,h):=\max_{v,u\in V(\Gamma)}\,|x_v-x_u|$;
        \item The curve $(\Gamma,h)$ is called {\em vertically stretched} (or {\em stretched} for brevity) if $|y_v-y_u|\gg \w(\Gamma,h)$ for any pair of vertices $v$ and $u$ lying on different floors of $(\Gamma,h)$.
    \end{enumerate}
\end{defn}

Simple stretched floor decomposed ({\em ssfd}) tropical curves arise naturally when imposing vertically stretched point constraints. More precisely, set $n:=|\nabla|+g-1$, and consider $n$ tropically general points $q_i=(x_i,y_i)\in\RR^2$ in a vertically stretched position, i.e., $|y_i-y_j|\gg\max_{k,s}\{x_k-x_s\}$ for all $i\ne j$. Then, by \cite[\S 5.1]{BM08}, any stable parametrized tropical curve $h\:\Gamma\to\RR^2$ of degree $\nabla$ such that $\{q_i\}\subset h(\Gamma)$ is simple, floor decomposed, and stretched. Moreover, the image of any elevator and any floor of $\Gamma$ contains a unique point $q_i$.

\begin{rem}
    Unlike in the case of general floor decomposed curves, the floors of stretched floor decomposed curves are naturally ordered by their height from the bottom to the top; cf.\ Figure~\ref{fig:comtype}. 
\end{rem}

\subsection{Moduli of parametrized tropical curves}\label{subsubsec:moduli of parametrized tropical curves}

Denote by $M_{g,n, \nabla}^{\trop}$ the moduli space of stable para\-metrized tropical curves of genus $g$ and degree $\nabla$, which have exactly $n$ contracted legs. If $n>0$, we will always assume that the first $n$ legs $l_1,\dotsc, l_n$ are the contracted ones. According to \cite{ACP}, $M_{g,n, \nabla}^{\trop}$ is a generalized polyhedral cone complex with integral structure. Recall that a {\em combinatorial type} $\Theta$ of a parametrized tropical curve $(\Gamma,h)$ is the following datum: the underlying weighted graph $\GG$ with ordered legs, equipped with the collection of slopes $\frac{\partial h}{\partial \vec e}$ for $e \in \overline E(\GG)$. The relative interiors of polyhedral cones in $M_{g,n, \nabla}^{\trop}$ are given by the subsets $M_{[\Theta]} = M_{\Theta}/\mathrm{Aut}(\Theta)$ for different combinatorial types $\Theta$. Here $M_{\Theta}$ is the interior of a polyhedron  $\overline M_\Theta$ in $\RR^{E(\GG)} \times N_\RR^{V(\GG)}$ parametrizing tropical curves $h\: \Gamma \to N_{\RR}$ of type $\Theta$: a curve $(\Gamma,h)$ corresponds to the point $\left((\ell(e))_e,(h(v))_v\right)\in M_\Theta\subset \RR^{E(\GG)} \times N_\RR^{V(\GG)}$; see, e.g., \cite[\S~3.1.4]{CHT23} for more details. We call the subsets $M_{[\Theta]} \subset M_{g,n, \nabla}^{\trop}$ the {\em strata} of the tropical moduli space. Set $M_{g,\nabla}^{\trop}:=M_{g,0,\nabla}^{\trop}$. There exists a natural forgetful map $M_{g,n,\nabla}^\trop\to M_{g,\nabla}^\trop$ given by removing the contracted legs and then stabilizing.

\begin{defn} \label{def:nice}
    A stratum $M_{[\Theta]}\subseteq M_{g,n, \nabla}^{\trop}$ is {\em nice} if it is regular (that is, of the expected dimension), and the underlying graph of $\Theta$ is weightless and 3-valent; $M_{[\Theta]}$ is {\em a simple wall} if it is regular and the underlying graph is weightless and 3-valent except for a unique 4-valent vertex. A stratum $M_{[\Theta]}\subset M_{g,\nabla}^\trop$ is called {\em simple stretched floor decomposed (ssfd)} if it contains an ssfd curve. By a slight abuse of language, we sometimes refer to $M_\Theta$ as a stratum, and say that $M_\Theta$ is nice, a simple wall, or ssfd if $M_{[\Theta]}$ is.
\end{defn}

By \cite[Proposition 2.23]{Mik05}, simple parametrized tropical curves are contained in nice strata; that is, their combinatorial type is regular. In particular, ssfd strata are nice. Each simple wall $M_{[\Theta]}$ is contained in the closure of exactly 3 strata, all of them nice, corresponding to the three ways of splitting the 4-valent vertex into a pair of 3-valent vertices. Finally, note that the forgetful map takes a nice stratum to a nice stratum, and takes a simple wall to a nice stratum if the 4-valent vertex is adjacent to a tree with contracted leaf-ends, and to a simple wall otherwise.

\subsection{Families of parametrized tropical curves}\label{subsubsec:family of parametrized tropical curves} 

We sketch the terminology and refer the reader to \cite[\S~2]{CHT24a} for formal definitions.
Let $\Lambda$ be a polyhedral complex with an integral affine structure (e.g., a tropical curve). By a {\em family} of parametrized tropical curves over $\Lambda$, we mean a continuous family of curves $h \colon \Gamma_\Lambda \to N_\RR$ with fibers $h_\kappa \colon \Gamma_\kappa \to N_\RR$ for $\kappa\in\Lambda$ such that the following hold:
\begin{enumerate}
    \item The degree and the number of contracted legs are the same for all fibers. 
    \item Along the interior of every face of $\Lambda$, the combinatorial type of the fibers is constant.
    \item For any two faces $W,W'$ of $\Lambda$ such that $W'$ is a face of $W$, we are given a contraction of weighted graphs with marked legs $\varphi_{W,W'}\colon \GG_{W} \to \GG_{W'}$ of the underlying graphs of the fibers. 
    \item The length of an edge $\gamma \in E(\GG_W)$ in $\Gamma_\kappa$ along a face $W$ of $\Lambda$ is given by an integral affine function. Here we identify edges over the interior of $W$ with those over a face $W'$ of $W$ via $\varphi_{W,W'}$, setting their length to zero in case they get contracted.
    \item Similarly, the coordinates of the images $h(u)$ for a vertex $u$ of $\GG_W$ in $h_\kappa(\Gamma_\kappa)$ along a face $W$ of $\Lambda$ are given by integral affine functions, where we again identify vertices over the interior of $W$ with those over faces $W'$ of $W$ via $\varphi_{W,W'}$.
\end{enumerate}
Any family of parametrized tropical curves $h\:\Gamma_\Lambda \to N_\RR$ induces a piecewise integral affine map \[\alpha\: \Lambda \to M_{g, n, \nabla}^{\trop}\] by sending $\kappa \in \Lambda$ to the point parametrizing the isomorphism class of the stabilization of the fiber $h_\kappa\: \Gamma_\kappa \to N_\RR$. Furthermore, if $W^\circ$ is the interior of a face $W$ and $\alpha(W^\circ)\subseteq M_{[\Theta]}$, then $\alpha|_{W^\circ}$ factors through $M_\Theta$.

\subsection{Tropicalization of a parametrized curve}\label{subsec:tropicalization of parametrized tropical curves}

There are several tropicalization constructions in the literature. The one we use goes back to \cite{Tyo12}; see also \cite[\S~4.2]{CHT23}. Let $(C, \sigma_\bullet,f\: C \to S)$ be a parametrized curve and $(C^0,\sigma_\bullet^0) \to \Spec(K^0)$ a prestable model of $(C, \sigma_\bullet)$. Denote by $\widetilde C$ the fiber of $C^0$ over the closed point of $\Spec(K^0)$. The {\em tropicalization} $\trop(C)$ of $(C, \sigma_\bullet)$ with respect to the given model is the following tropical curve $\Gamma=(\GG,\ell,g)$: the graph $\GG$ is the dual graph of the reduction $(\widetilde C,\widetilde{\sigma}_\bullet)$, i.e., the vertices of $\GG$ correspond to irreducible components of $\widetilde C$, the edges to nodes, the legs to marked points, and the natural incidence relation holds; the weight $g(v)$ is the geometric genus of the corresponding component of the reduction $\widetilde C$; and the length function is given by $\ell(e):=\val(\lambda)$, where $\lambda\in K^{00}$ is such that the total space of $C^0$ is given by $xy=\lambda$ \'etale locally at the node corresponding to the edge $e$. Although $\lambda$ depends on the \'etale neighborhood, its valuation does not, and hence $\ell(e)$ is well-defined. Finally, the order of the marked points induces the order on the legs of $\GG$.

To define the parametrization $h\:\trop(C)\to N_\RR$, it suffices to specify its values at the vertices and its slopes along the legs. Let $\widetilde C_v$ be the irreducible component of $\widetilde C$ corresponding to a vertex $v$. For any $m \in M$, consider the monomial function $x^m$ and its pullback to $C^0$. Since $f^{-1}(T)$ is dense in $C$, the pullback $f^*(x^m)$
is a non-zero rational function on $C^0$, and therefore, there exists $\lambda_m \in K^\times$ such that $\lambda_mf^*(x^m)$ is an invertible function at the generic point of $\widetilde C_v$. Set $h(v)(m):=\val(\lambda_m)$. Since $\lambda_m$ is unique up to an invertible element of $K^0$, the value of $h(v)(m)$ is well-defined. Plainly, $h(v)$ is linear in $m$. Thus, $h(v)\in N_\RR$. The slope of $h$ along the leg $l$ corresponding to a marked point $\sigma_l$ is defined similarly by setting $\frac{\partial h}{\partial \vec l}(m):=-\mathrm{ord}_{\sigma_l} f^*(x^m)$ for all $m\in M$. Plainly, the slope defined this way is an element of $N$, and it is not difficult to verify that the balancing condition holds; see, e.g., \cite[Lemma~2.23]{Tyo12}. 

The parametrized tropical curve $(\trop(C),h)$ constructed this way is called the {\em tropicalization} of $(C, \sigma_\bullet,f\: C \to S)$ with respect to the model $(C^0,\sigma_\bullet^0)$. The tropical curve $\trop(C)$ is independent of the parametrization $f$ and depends only on the model. If the model is the stable model of $(C,\sigma_\bullet)$, then the corresponding tropicalization is called simply {\em the tropicalization} of $C$ (resp.\ of $f\: C\to S$). There is a natural {\em tropicalization map} \[\trop\: C(K)\setminus \bigcup_i \sigma_i \to \trop(C)\] from the set of non-marked $K$-points of $C$ to $\trop(C)$. For a $K$-point $\eta\in C(K)\setminus \bigcup_i \sigma_i$, its image $\trop(\eta)$ is defined as follows. Add $\eta$ to the sequence of marked points, and consider the minimal modification $(C^0)'\to C^0$ such that $(C^0)'$ is a prestable model of $(C,\sigma_\bullet,\eta)$. Let $l$ be the leg of the new tropicalization $\trop(C)'$ corresponding to the marked point $\eta$. Then, $\trop(C)'$ is obtained from $\trop(C)$ by attaching $l$ either to an existing vertex $w$ or to a new 2-valent vertex $w$ of weight zero that splits an edge or a leg of $\trop(C)$. The tropicalization map $\trop$ sends $\eta$ to $w$. 

\begin{rem}
    Any integral curve $C \subset S$ disjoint from the zero-dimensional orbits in $S$ defines a parametrized curve $f \colon C^\nu \to S$, where $C^\nu$ denotes the normalization of $C$ and the marked points are the preimages of the boundary divisor $\partial S$ with some choice of labeling. Consequently, we define the {\em tropicalization} of $C$ to be the tropicalization of $f \colon C^\nu \to S$, which is uniquely defined by $C$ up to a change of the ordering of the legs. By abuse of notation, we denote this tropicalization by $\trop(C)$. 
\end{rem}

\subsection{Tropicalization in families} \label{subsec:tropicalization in families}

The tropicalization construction of the previous section can be performed in families after an appropriate modification of the base of the family. For families over one-dimensional bases, this was carried out in \cite{CHT23}, and in general in \cite{CHT24a}. While we do use results of \cite{CHT24a} that rely on the higher-dimensional construction, in the current paper we will make explicit use only of tropicalizations over one-dimensional bases. Thus, we sketch here only this case. 

Let $(\cC \to B, \sigma_{\bullet}, f\colon \cC \rightarrow S)$ be a family of parametrized curves over an integral curve $B$ as in \S~\ref{subsec: parametrized curves}. In order to construct the tropicalization of $\cC \to B$, we will need to pull back the family along several alterations of the base $B$. To ease the notation, we will keep calling the base $B$ even after these alterations. See Step 2 in the proof of \cite[Theorem~5.1]{CHT23} for more details.

First, we pass to the normalization of $B$ so that we may assume that $B$ is smooth. By simultaneous stable reduction, up to an alteration of the base, we may then assume that the family $\mathcal C \to B$ extends to a split stable model $\mathcal C^0 \to B^0$ of marked curves over a projective, (pre)stable model $B^0$ of $B$ over $K^0$. We replace $B \subset B^0$ by the generic fiber of $B^0$ over $K$ and let $B' \subset B$ be the locus of $K$-points over which $f$ extends to a family of parametrized curves. Since $B$ is one-dimensional, the complement of $B'$ is a finite collection of points $\tau_\bullet$, which we may view as marked points of $B$. Up to another alteration, we may assume that $B^0$ is a (pre)stable model for $B$ also after adding the new marked points.

The tropicalization $\Lambda$ of the curve $B$ with respect to the model $B^0$ is the base of the tropicalization of $(\cC \to B, \sigma_{\bullet}, f\colon \cC \rightarrow S)$. It remains to construct the family of parametrized tropical curves over $\Lambda$. For a point $\eta \in B'(K)$, set $q:=\trop(\eta)\in\Lambda$. By \cite[Theorem~4.6]{CHT23}, the tropicalization of $f|_{\cC_\eta}\colon \cC_{\eta} \to S$ depends only on $q$ and the pointwise tropicalizations $h_q\: \Gamma_q \to N_\RR$ glue to a family of parametrized tropical curves. Using the splitness of the model $\cC^0 \to B^0$, one verifies the five axioms of \S~\ref{subsubsec:family of parametrized tropical curves} for the fibers over the rational points of $\Lambda$, and then uniquely extends the family by continuity to the rest of $\Lambda$.

\begin{defn}
    We call the family of parametrized tropical curves $h \colon \Gamma_\Lambda \rightarrow N_\RR$ obtained as above a {\em tropicalization} of the family of parametrized curves $f \colon \mathcal C \to S$.
\end{defn}

\subsection{A tropicalization of the Severi variety and its components}\label{subsec:tropseveri}
\begin{nota}\label{nt:set of trops} 
    Let $\Delta\subset M_\RR$ be a lattice polygon, $V\subseteq V^\irr_{g,\Delta}$ an irreducible component of the Severi variety, and $\nabla$ the reduced tropical degree dual to $\Delta$. We denote by $\Sigma_V\subseteq M_{g,\nabla}^\trop$ the collection of all tropicalizations of parametrized curves $(C^\nu,\sigma_\bullet, f\:C^\nu\to S_\Delta)$, where $[C]\in V(K)$ is a curve intersecting the boundary divisor $\partial S_\Delta$ transversely, $C^\nu$ its normalization, $f\:C^\nu\to S_\Delta$ the natural map, and $\sigma_\bullet\subset C^\nu$ a marking of the preimage of the boundary divisor.
\end{nota}

One can think about $\overline{\Sigma}_V$ as a tropicalization of the component $V$. The following properties of $\overline{\Sigma}_V$ follow easily from \cite{CHT24a}.

\begin{lem}\label{lem:surjectivity into nice cones}
    Suppose $\Delta\subset M_\RR$ is a lattice polygon. Then,
    \begin{enumerate}
        \item $\Sigma_V$ is the set of rational points of a rational polyhedral subcomplex in $M^{\trop}_{g,\nabla}$;
        \item If $M_{[\Theta]}$ is nice and $\dim(\Sigma_V \cap M_{[\Theta]})=\dim M_{[\Theta]}$, then $M_{[\Theta]}\subseteq \overline{\Sigma}_V$;
        \item If $M_{[\Theta]}$ is a simple wall and $\overline{\Sigma}_V$ contains one of its adjacent nice strata, then $\overline{\Sigma}_V$ contains the other two as well. 
    \end{enumerate} 
\end{lem}
    
\begin{proof}
    To apply the results of \cite{CHT24a}, we must first replace the family of curves over $V$ with a family of parametrized curves while keeping the set $\Sigma_V$ unchanged.

    Let $V'\subset V$ be the locus of curves intersecting the boundary divisor transversely, and $\mathcal C \to V'$ the universal family of curves over $V'$. According to \cite[\S~5.10]{dJ96}, up to an alteration of $V'$, we may assume that $\mathcal C$ is generically equinormalizable. Hence, the geometric genus of the generic fiber of $\mathcal C$ agrees with the geometric genus of the fibers over the closed points. Pulling $\mathcal C$ back to the normalization of $V'$ if necessary, we may assume that $V'$ is normal. Consequently, the restriction of $\mathcal C$ to the localization of $V'$ at any closed point is equinormalizable by \cite[Theorem 4.2]{CL06}. Hence, the whole family $\mathcal C$ is equinormalizable. The normalization $\mathcal C^\nu \to V'$ of the total space of the family $\mathcal C \to V'$ then gives a family of smooth curves along with a map $\mathcal C^\nu\rightarrow S_\Delta$. It remains to mark the preimage of the boundary divisor, which can be done after yet another alteration of the base (cf.\ \cite[\S~5]{dJ96}). Note that the above process does not change $\Sigma_V$. Therefore, the lemma now follows from \cite[Corollary 3.2]{CHT24a}.
\end{proof}

\section{A lifting result}\label{sec:lifting}
The goal of the current section is to prove the following theorem, which provides a correspondence between the strata in the tropical moduli spaces containing simple curves of multiplicity one and irreducible components of the Severi variety. The uniqueness part of the theorem plays a central role in the proof of Theorem~B.

\begin{thm}\label{thm:existence and uniqueness}
    Let $g\ge 0$ be an integer, $\Delta\subset M_\RR$ a lattice polygon, and $\nabla$ the reduced tropical degree dual to $\Delta$. Let $M_{[\Theta]}\subset M_{g,\nabla}^\trop$ be a nice stratum containing the isomorphism class of a simple parametrized tropical curve $(\Gamma,h)$ of multiplicity one. Then there exists a unique irreducible component $V\subseteq V^\irr_{g,\Delta}$ such that $M_{[\Theta]}\subset \overline \Sigma_V$.
\end{thm}

As explained in the introduction, it is not surprising that the multiplicity-one strata satisfy the unique lifting property. Our proof proceeds as follows. After picking $|\nabla|+g-1$ tropically general points on $h(\Gamma)$ together with a general lift of these points to the algebraic torus, we apply an algebraic version of Mikhalkin's Correspondence Theorem that works in any characteristic. Since a multiplicity-one curve satisfying a tropical constraint admits a unique lift satisfying the given lift of the constraint, this allows us to conclude not only the existence of the desired component but also its uniqueness, thanks to the following proposition. 

\begin{prop}\label{prop:lifting genral curve through points to all comp}
    Let $g\ge 0$ be an integer, $\Delta\subset M_\RR$ a lattice polygon, and $\nabla$ the reduced tropical degree dual to $\Delta$. For $n=|\nabla|+g-1$, let $q_1,\dotsc, q_n\in N_\QQ$ be points in tropically general position and $(\Gamma,h)$ a parametrized tropical curve in $M_{g,\nabla}^\trop$ passing through the $q_i$'s, i.e., $\{q_i\}_{i=1}^n\subset h(\Gamma)$. Then, there exist $p_i\in T(K)$ with $\trop(p_i)=q_i$ such that for \emph{any} irreducible component $V\subseteq V^\irr_{g,\Delta}$ that contains a curve tropicalizing to $(\Gamma,h)$, there exists a curve $[C]\in V$ such that $p_i\in C$ for all $i$ and $\trop(C)=(\Gamma,h)$ for a suitable choice of marking of the preimage of the boundary divisor $\partial S_\Delta$ on the normalization $C^\nu$.
\end{prop}

\subsection{A technical Lemma}\label{subsec:technical preparations}

The main ingredient needed for the proof of Proposition~\ref{prop:lifting genral curve through points to all comp} is a technical lemma, which, roughly speaking, says the following. Let $h \colon \Gamma_\Lambda \to N_\RR$ be the tropicalization of a family of parametrized curves $f \colon \cC \to S$ over a base $B$, and $w$ a vertex of $\Lambda=\trop(B)$. Assume that the fiber $(\Gamma_w,h_w)$ of the tropical family passes through a point $q \in N_\QQ:=N\otimes_\ZZ\QQ$, i.e., $q\in h_w(\Gamma_w)$, and that locally at $q$ the family sweeps out a two-dimensional region of $N_\RR$. Then for almost any lift $p\in T\subset S$ of $q$, one can find $b\in B$ such that $f \colon \cC_b \to S$ passes through $p$ and has tropicalization $(\Gamma_w, h_w)$. That is, not only can we lift $(\Gamma_w, h_w)$ to an element of the algebraic family, but we can lift it in such a way that the image of the algebraic curve in $S$ passes through almost any given lift of $q$. 

Let us now explain what we mean by ``almost any lift''. The locus of points in the algebraic torus with a given tropicalization is naturally a torsor under the action of $T(K^0)$. Picking a point $t$ with $\trop(t)=q$ trivializes this torsor and induces a scaled reduction map to the torus over the residue field $T(\widetilde{K})$ given by $p\mapsto \widetilde{t^{-1}p}$. Thus, by almost any lift $p$ of $q$, we simply mean that $p$ belongs to the preimage of an open dense subset under the scaled reduction map to $T(\widetilde{K})$.

To simplify the notation, we will state the result only in the case of a one-dimensional base, which is sufficient for our purposes.

\begin{lem}\label{lem:lifting} 
Let $S$ be a toric surface, $T\subset S$ the dense open orbit, and $f \colon \cC \to S$ a family of parametrized curves with marked points $\sigma_\bullet$ over a smooth proper base curve $(B,\tau_\bullet)$ defined over $B\setminus \tau_\bullet$. Let $h\: \Gamma_\Lambda\to N_\RR$ be its tropicalization, $\alpha\:\Lambda\to M^{\trop}_{g,n,\nabla}$ the induced map, and $\vec e\in \overline{E}(\Lambda)$. Let $\Theta$ be a combinatorial type for which $\alpha(e^\circ)\subseteq M_{[\Theta]}$ and $\vec \gamma\in\overline{E}(\GG_\Theta)$ an edge or a leg. Set $w:=\ft(\vec e)$ and $u:=\ft(\vec\gamma)$, and pick $t\in T(K)$ such that $\trop(t)=h_w(u)$. Assume that:
\begin{enumerate}
    \item $\alpha(w)\in M_{[\Theta]}$,
    \item the partial derivatives $\frac{\partial h_\bullet(u)}{\partial \vec e}$ and $\frac{\partial h_w(\bullet)}{\partial \vec\gamma}$ are linearly independent. \label{liftcond2}
\end{enumerate}
Then, there exists an open dense $U\subseteq T(\widetilde{K})$ such that for any $p\in T(K)$ satisfying $\trop(p)=h_w(u)$ and $\widetilde{t^{-1}p}\in U$, there exists $b\in B(K)$ such that $\trop(b)=w$ and $p\in f(\cC_b)$.
\end{lem}

\begin{rem}
    The linear independence condition implies that the family sweeps out a two-dime\-nsional region of $N_\RR$ containing $q:=h_w(u)$. At first glance, the linear independence condition may look stronger, but in fact, the two conditions are equivalent if the vertices $w$ and $u$ are chosen wisely.
\end{rem}

\begin{figure}[ht]  
    \begin{tikzpicture}[x=0.8pt,y=0.8pt,yscale=-0.6,xscale=0.6]
        \import{./}{figlift2.tex}
    \end{tikzpicture}   
    \caption{A geometric illustration of the linear independence assumption of Lemma~\ref{lem:lifting}. The image of the curves in the family deforms from $h_w(\Gamma_w)$ to $h_v(\Gamma_v)$ and sweeps out a two-dimensional region in a neighborhood of $h_w(u)$.}
    \label{figlift2}
\end{figure}

\begin{proof}
    In the proof, we will assume that $\vec \gamma\in E(\GG_\Theta)$ is an edge. The case of $\vec\gamma$ being a leg is similar and is left to the reader. After rescaling the parametrization map $f$ by $t^{-1}$, we may assume that $h_w(u)=0$ and $t=1$. Let $(\cC^0 \to B^0,\sigma_\bullet ^0)$ be a split stable model with respect to which the tropicalization is constructed, see \S~\ref{subsec:tropicalization in families}. As usual, to construct such a model, one may need to replace $B$ with an alteration first. After modifying the model $B^0$, we may assume that $\vec e\in E(\Lambda)$ is an edge, and that $\alpha(\fh(\vec e))\in M_{[\Theta]}$. Set $v:=\fh(\vec e)$, and let $\widetilde{B}_{w}$ and $\widetilde{B}_{v}$ be the components of the reduction of $B$ corresponding to $w$ and $v$, respectively. Let $\widetilde{B}_{\vec e}\in \widetilde{B}_{w}\cap \widetilde{B}_{v}$ be the node of the reduction corresponding to $\vec e$. 
    
    Set $u':=\fh(\vec\gamma)$. Denote by $\widetilde{\cC}_{w,u}$ the component of $\widetilde{\cC}$ over $\widetilde{B}_{w}$ corresponding to the vertex $u$, and similarly for $\widetilde{\cC}_{w,u'}$, $\widetilde{\cC}_{v,u}$, $\widetilde{\cC}_{v,u'}$, $\widetilde{\cC}_{\vec e,u}$, and $\widetilde{\cC}_{\vec e,u'}$. Finally, let $\rho_{\vec\gamma}\:\widetilde{B}_{w}\cup \widetilde{B}_{v}\to \widetilde{\cC}$ be the section, whose image is the node corresponding to the edge $\vec\gamma$; see Figure~\ref{figlift1} for an illustration. Such a section exists since the model is split. Set $s:=\rho_{\vec\gamma}(\widetilde{B}_{\vec e})$.

    \begin{figure}[ht]  
        \begin{tikzpicture}[x=0.8pt,y=0.8pt,yscale=-0.6,xscale=0.6]
            \import{./}{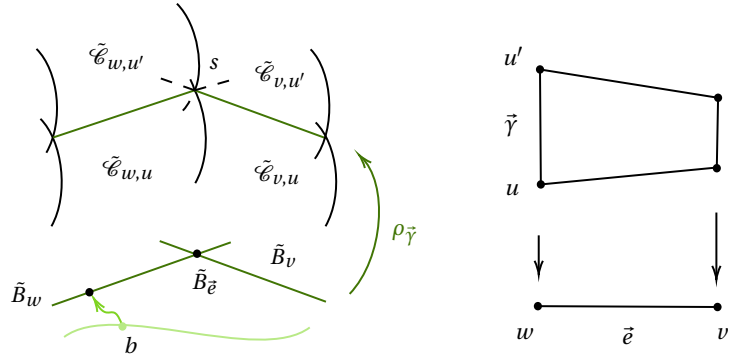}
        \end{tikzpicture}   
        \caption{The setting in the proof of Lemma~\ref{lem:lifting}. On the left the local picture close to $\widetilde B_{\vec e}$ in $B^0$. On the right, the local picture of its tropicalization.}
        \label{figlift1}
    \end{figure}

    Step 1: {\em The construction of $U$.} By our assumption, $h_w(u)=0$. Thus, the monomial functions induce a rational map $\tilde{f}_{w,u}\:\widetilde{\cC}_{w,u}\dashrightarrow T(\widetilde{K})$. Let us first show that this map is dominant. Pick a pair of monomials $x^{m_1},x^{m_2}$ such that 
    \begin{equation}\label{eq:dualsemibasis}
        \left(m_1,\frac{\partial h_w(\bullet)}{\partial \vec\gamma}\right)<0, \left(m_2,\frac{\partial h_\bullet(u)}{\partial \vec e}\right)<0,\left(m_2,\frac{\partial h_w(\bullet)}{\partial \vec\gamma}\right)=\left(m_1,\frac{\partial h_\bullet(u)}{\partial \vec e}\right)=0.
    \end{equation} 
    Such a pair exists since $\frac{\partial h_\bullet(u)}{\partial \vec e}$ and $\frac{\partial h_w(\bullet)}{\partial \vec\gamma}$ are linearly independent. It is sufficient to show that the map $(x^{m_1},x^{m_2})\:\widetilde{\cC}_{w,u}\dashrightarrow (\widetilde{K}^\times)^2$ is dominant.

    The functions $x^{m_1}$ and $x^{m_2}$ are regular at the generic points of $\widetilde{\cC}_{w,u}$, $\widetilde{\cC}_{w,u'}$, $\widetilde{\cC}_{v,u}$, and $\widetilde{\cC}_{v,u'}$. Indeed, recall that we are working over the algebraic closure $K$ of a complete discretely valued field $F$. Therefore, after replacing $F$ with a finite extension, we may assume that the family $\cC$ is defined over $F$, whose uniformizer has valuation one. By the construction of the tropicalization of curves, for each $m_i$, the order of vanishing of $x^{m_i}$ along $\widetilde{\cC}_{w,u}$ is 
    $$\mathrm{ord}_{\widetilde{\cC}_{w,u}}(x^{m_i})=-(m_i,h_w(u))=0,$$ 
    and similarly for the three other components. Thus, by \eqref{eq:dualsemibasis}, $\mathrm{ord}_{\widetilde{\cC}_{w,u'}}(x^{m_1})>0$, $\mathrm{ord}_{\widetilde{\cC}_{w,u'}}(x^{m_2})=0$, $\mathrm{ord}_{\widetilde{\cC}_{v,u}}(x^{m_1})=0$, $\mathrm{ord}_{\widetilde{\cC}_{v,u}}(x^{m_2})>0$, $\mathrm{ord}_{\widetilde{\cC}_{v,u'}}(x^{m_1})>0$, and $\mathrm{ord}_{\widetilde{\cC}_{v,u'}}(x^{m_2})>0$. As a result, $x^{m_1}$ and $x^{m_2}$ are regular in a neighborhood $\Omega$ of the point $s$ by the algebraic Hartogs lemma. Furthermore, the function $x^{m_1}$ vanishes along $\Omega\cap\rho_{\vec\gamma}(\widetilde{B}_{w})$ and is invertible at the generic point of $\Omega\cap\widetilde{\cC}_{\vec e,u}$, while the function $x^{m_2}$ vanishes along $\Omega\cap\widetilde{\cC}_{\vec e,u}$ and is invertible at the generic point of $\Omega\cap\rho_{\vec\gamma}(\widetilde{B}_{w})$. Therefore, the morphisms $x^{m_1}\:\Omega\cap\widetilde{\cC}_{\vec e,u}\to \AAA^1$ and $x^{m_2}\:\Omega\cap\rho_{\vec\gamma}(\widetilde{B}_{w})\to \AAA^1$ are dominant. 

    Consider the morphism $(x^{m_1},x^{m_2})\:\Omega\cap\widetilde{\cC}_{w,u}\to \AAA^2$. It is sufficient to show that this morphism is dominant. Assume to the contrary that it is not dominant. Then its image belongs to an irreducible curve $Z\subset \AAA^2$. Since $x^{m_2}$ is invertible at the generic point of $\Omega\cap\rho_{\vec\gamma}(\widetilde{B}_{w})$ and $x^{m_1}$ vanishes at it, the curve $Z$ must contain the generic point of the $y$-axis in $\AAA^2$. In a similar manner, it necessarily contains the generic point of the $x$-axis, which is absurd since $Z$ is an irreducible curve. Therefore, $(x^{m_1},x^{m_2})\:\Omega\cap\widetilde{\cC}_{w,u}\to \AAA^2$ is dominant, and hence so are $(x^{m_1},x^{m_2})\:\widetilde{\cC}_{w,u}\dashrightarrow (\widetilde{K}^\times)^2$ and $\tilde{f}_{w,u}\:\widetilde{\cC}_{w,u}\dashrightarrow T(\widetilde{K})$. 

    By dimension theory, since $\tilde{f}_{w,u}\:\widetilde{\cC}_{w,u}\dashrightarrow T(\widetilde{K})$ is a dominant rational map between algebraic surfaces, there exists an open dense subset $U\subseteq T(\widetilde{K})$ such that all fibers over $U$ contain isolated points in $\widetilde{\cC}_{w,u}^\circ$, where $\widetilde{\cC}_{w,u}^\circ$ denotes the complement of the intersection of $\widetilde{\cC}_{w,u}$ with other components of the reduction and with the reductions of the marked points $\sigma_\bullet$ that are not mapped to the dense orbit $T(K)$. 

    Step 2: {\em The set $U$ satisfies the assertion of the lemma}. Let $p\in T(K)$ be such that $\trop(p)=0$. Then $p$ admits a unique integral model $p^0\in T(K^0)$. Assume that its reduction $\widetilde{p}$ belongs to $U$. We must show that there exist $b\in B(K)$ and $p'\in \cC_b(K)$ such that $\trop(b)=w$ and $f(p')=p$. Let us first explain briefly the idea of the argument, as it is a bit technical. We will look for $p'$ among the points whose reduction belongs to $\widetilde{\cC}_{w,u}^\circ$, and define $b$ to be its projection to the base $B$. The surjectivity of $\tilde{f}_{w,u}\:\widetilde{\cC}_{w,u}\dashrightarrow T(\widetilde{K})$ over $U$ will allow us to construct the desired reduction, and then we will use a Hensel's Lemma type argument to show that it extends to a $K^0$-point, whose base change to $K$ is mapped to $p$ by the map $f$. Here is the precise argument. 
    
    After replacing $F$ with a finite extension, we may assume that the models of $(\cC^0\to B^0, \sigma_\bullet^0)$ and the point $p^0$ are defined over $F^0$. Let us construct $b$ and $p'$ as asserted. Consider the open subset $\Omega':=(\cC\setminus\sigma_\bullet)\cup \widetilde{\cC}_{w,u}^\circ\subset\cC^0$ and the morphism $f\:\Omega'\to T$. Then, by the construction of $U$, the pullback of $\tilde{p}\in T(\widetilde{F})$ has an isolated point $z\in \widetilde{\cC}_{w,u}^\circ\subset\Omega'$. Thus, $f$ is flat at $z$ by \cite[Theorem~23.1]{Mat89}, and hence so is the base change $f^{-1}(p^0)\to p^0=\Spec(F^0)$. Therefore, by a result of Mumford (cf.\ \cite[Proposition~14.5.10]{grothendieck1966elements}) there exists a finite extension $F\subseteq L\subset K$ and a section $\Spec(L^0)\to f^{-1}(p^0)\times_{p^0} \Spec(L^0)$ mapping the closed point to $z$. Let $p'\in \cC$ be the image of $\Spec(L)$, and $b$ be its projection to the base curve $B$. Then $p'\in \cC_b$ and $f(p')=p\in T$. Furthermore, since the reduction of $p'$ belongs to $\widetilde{\cC}_{w,u}^\circ$, the reduction of $b$ belongs to $\widetilde{B}_w$, and therefore, $\trop(b)=w$, which completes the proof.
\end{proof}

\begin{rem}
    If $g(w)=g(u)=0$, the vertex $w$ is 2-valent, and $u$ is 3-valent and is adjacent to a contracted leg, then one can choose $U=T(\widetilde{K})$, and thus, for any $p\in \trop^{-1}(h_w(u))$ there exists $b\in B(K)$ such that $p\in f(\cC_b)$. To see this, one notices that under such assumptions, $\widetilde{\cC}_{w,u}^\circ$ is isomorphic to the two-dimensional torus. Since any invertible function on an algebraic torus is a scalar multiple of a monomial, it follows that up to a translation by an element of $T(\widetilde{K})$, the morphism $\tilde{f}_{w,u}\:\widetilde{\cC}_{w,u}^\circ\to T(\widetilde{K})$ is a toric map, i.e., corresponds to a linear map between lattices of cocharacters. Since $\tilde{f}_{w,u}$ is dominant, the corresponding linear map is surjective, and therefore, so is $\tilde{f}_{w,u}$. Furthermore, since both tori are of dimension two, the map $\tilde{f}_{w,u}$ is finite, and hence, one can choose $U$ to be the whole of $T(\widetilde{K})$.
\end{rem}

\subsection{Proof of Proposition~\ref{prop:lifting genral curve through points to all comp}}
    Let $0\le k\le n$ be the maximal integer for which there exist $p_i\in T$ such that $\trop(p_i)=q_i$ and for any irreducible component $V\subseteq V^\irr_{g,\Delta}$ that contains a curve tropicalizing to $(\Gamma,h)$, there exists a curve $[C]\in V$ such that $p_i\in C$ for all $1\le i\le k$ and $\trop(C)=(\Gamma, h)$ for a suitable choice of marked points on the normalization $C^\nu$ mapping to the boundary divisor $\partial S_\Delta$. Assume to the contrary that $k<n$. Let us fix a corresponding tuple of points $\{p_i\}_{i=1}^k$, and pick $t\in T(K)$ such that $\trop(t)=q_{k+1}$.
    
    Let $V$ and $[C]\in V$ be as above. For all $i>k$, pick any $p_i\in C$ satisfying $\trop(p_i)=q_i$. Such points exist since $q_i\in N_\QQ$. Consider the locus $B\subset V$ of curves passing through the points $\{p_i\}_{i\ne k+1}$. Since $n=|\nabla|+g-1=\dim(V)$ and $\{q_i\}_{i=1}^n$ are in tropically general position, the locus $B$ is a curve, as otherwise $\dim(B)\ge 2$, and therefore $\trop(C_b)$ would pass through $n+1$ points in tropically general position for some $b\in B$, which is absurd. 
    
    We will employ Lemma~\ref{lem:lifting} to get the desired contradiction. To do so, we shall first adjust the setting accordingly. After replacing $B$ with its smooth projective model and marking finitely many points $\tau_\bullet$ on it, we may assume that $B$ is a smooth projective curve, and that the tautological family of curves over it is defined over the complement of the marked points $\tau_\bullet$. Furthermore, by passing to a finite covering, we may assume that the tautological family is equinormalizable over this complement.
    
    Let $\cC\to B\setminus \tau_\bullet$ be the normalization of the tautological family, $f\:\cC\to S_\Delta$ the natural map, and let us mark a point $\tau_0\in B(K)$ such that the fiber over it is the normalization of the curve $C$. After replacing $(B,\tau_\bullet)$ with yet another finite covering, we may assume also that there exist $n+|\nabla|$ sections $\sigma_i\:B\dashrightarrow \cC$ such that $\bigcup_{i>n}\sigma_i=f^{-1}(\partial S_\Delta)$ is the pullback of the boundary divisor, $\trop\left(\cC_{\tau_0},(\sigma_i(\tau_0))_{i>n}\right)=(\Gamma,h)$, $f(\sigma_i)=p_i$ for $k+1\ne i\le n$, and $f(\sigma_{k+1}(\tau_0))=p_{k+1}$. To summarize, we have constructed a one-parameter family $f \colon \cC \to S_\Delta$ of stable parametrized curves with marked points $\sigma_\bullet$ over a smooth proper base curve $(B,\tau_\bullet)$ defined over $B\setminus \bigcup_{i>0}\tau_i$. 

    Let $h\: \Gamma_\Lambda\to N_\RR$ be the tropicalization of the family $f \colon \cC \to S_\Delta$ and $\alpha\:\Lambda\to M^{\trop}_{g,n,\nabla}$ the induced map. By construction, $\alpha(\trop(\tau_0))$ maps to $(\Gamma,h)$ under the forgetful map $M^{\trop}_{g,n,\nabla}\to  M^{\trop}_{g,\nabla}$. Denote by $\Theta$ the combinatorial type of $\alpha(\trop(\tau_0))$. Then $M_{[\Theta]}$ is a nice stratum since $n=|\nabla|+g-1$ and $\{q_i\}_{i=1}^n$ are in tropically general position. Moreover, by a similar argument as above, a general member of the family $h\: \Gamma_\Lambda\to N_\RR$ does not pass through $q_{k+1}$. Thus, there exists a vertex $w\in V(\Lambda)$ such that $\alpha(w)=\alpha(\trop(\tau_0))$ but $q_{k+1}\notin h_\kappa(\Gamma_\kappa)$ for $\kappa\in e^\circ$ for some $\vec e\in \Star(w)$. Since $M_{[\Theta]}$ is nice, $\alpha(e^\circ)\subset M_{[\Theta]}$. Let $\vec l_{k+1}\in L(\GG)$ be the leg corresponding to the marked point $\sigma_{k+1}$. Set $u:=\ft(\vec l_{k+1})$, and let $\vec\gamma\in\Star(u)$ be an edge or a leg not contracted by $h$. It follows now that the assumptions of Lemma~\ref{lem:lifting} are satisfied.
    
    By Lemma~\ref{lem:lifting}, there exists an open dense $U_V\subseteq T(\widetilde{K})$ such that for any $p\in T(K)$ satisfying $\trop(p)=h_w(u)$ and $\widetilde{t^{-1}p}\in U_V$, there exists $b\in B(K)$ such that $\trop(b)=w$ and $p\in f(\cC_b)$. Since there are only finitely many irreducible components $V$, the intersection of their respective open dense sets $U_V$ remains open and dense. Let $p\in T(K)$ be such that $\widetilde{t^{-1}p}$ belongs to this intersection. Then after replacing $p_{k+1}$ with $p$, we obtain a $(k+1)$-tuple $\{p_i\}_{i=1}^{k+1}$ satisfying the assertion of the first paragraph of the proof, which contradicts the maximality of $k$.
\qed

\subsection{Proof of Theorem~\ref{thm:existence and uniqueness}}
    The proof is based on the algebraic version of Mikhalkin's Correspondence Theorem due to the third author \cite{Tyo12}. For the convenience of the reader, we state the main result of {\em loc. cit.} using the notation of the current paper.
    
    Let $h'\:\Gamma'\to N_\RR$ be a parametrized tropical curve of genus $g$ and degree $\nabla$ and with $n$ contracted legs, where $n=|\nabla|+g-1$. As usual, we assume that the contracted legs are the first $n$ legs. Fix any orientation of the bounded edges of $\Gamma'$, and set $w_i:=\ft(l_i)$. Consider the two-term complex $L^\bullet(\Gamma')$:
    \begin{equation}\label{eq:two-term-complex}
        \zeta\: N^{V(\Gamma')}\times\ZZ^{E(\Gamma')}\to N^{E(\Gamma')}\times N^n 
    \end{equation}
    where the linear map $\zeta$ is given by $x_w\mapsto \left(\left(\epsilon(e,w)x_w\right),\left(\delta(w,w_i)x_w\right)\right)$ for any vertex $w\in V(\Gamma')$ and $y_e\mapsto \left(y_e\frac{\partial h'}{\partial \vec e}, 0\right)$ for any edge $e\in E(\Gamma')$ with the chosen orientation. Here, $\delta$ is the Kronecker delta, and $\epsilon(e,w)=1$ if $\ft(e)=w$, $\epsilon(e,w)=-1$ if $\fh(e)=w$, and $\epsilon(e,w)=0$ otherwise. Notice that changing the orientation of the edges gives rise to a quasi-isomorphic two-term complex. In particular, if $\zeta$ is given by a square matrix, then the absolute value of its determinant is independent of the choice of the orientation. Now we can state the particular case of \cite[Theorem~6.2]{Tyo12} that we need in the current paper.

    \begin{thm}\label{thm:corshort}
        Let $p_1,\dotsc, p_n\in T(K)$ be a tuple of points. Set $q_i:=\trop(p_i)$, and let $h'\:\Gamma'\to N_\RR$ and $\zeta$ be as above. Assume that:
        \begin{enumerate}
            \item The curve $\Gamma'$ is weightless and trivalent;
            \item The curve $(\Gamma',h')$ satisfies the tropical constraint $q_\bullet$, i.e., $h(l_i)=q_i$ for all $1\le i\le n$;
            \item The determinant of the map $\zeta$ is non-zero and prime to the characteristic of $\widetilde{K}$.
        \end{enumerate}
        Then, the number of parametrized marked curves $f\: C\to S_\Delta$ tropicalizing to $h'\:\Gamma'\to N_\RR$ and satisfying the constraint $p_\bullet$ is given by $|\det(\zeta)|$. In particular, if $|\det(\zeta)|=1$, then $h'\:\Gamma'\to N_\RR$ admits a unique algebraic lift. 
    \end{thm}
    
    \begin{rem}
        The assumption on the determinant of $\zeta$ implies assumptions 3, 4, 5, and 6 of \cite[Theorem~6.2]{Tyo12}, and therefore the main result of {\em loc. cit.} applies.
    \end{rem}
    
    Let us now return to Theorem~\ref{thm:existence and uniqueness} and prove it. Since $(\Gamma,h)$ is simple, the dimension of $M_{[\Theta]}$ is the expected one $\dim(M_{[\Theta]})=|\nabla|+g-1$. Thus, after perturbing $(\Gamma,h)$ slightly, we may assume that it passes through $n=|\nabla|+g-1$ points $q_1,\dotsc, q_n\in N_\QQ$ in tropically general position. In particular, $(\Gamma,h)\in M_\Theta$ corresponds to a general rational point. Pick algebraic lifts $p_i\in T(K)$ of the $q_i$'s satisfying the assertion of Proposition~\ref{prop:lifting genral curve through points to all comp}. Then the points $\{p_i\}$ are in general position since so are the $q_i$'s. Therefore, any curve in $V^\irr_{g,\Delta}$ that passes through the points $\{p_i\}_{i=1}^n$ belongs to a single irreducible component of the Severi variety. Moreover, if $f\:C\to S_\Delta$ tropicalizes to $(\Gamma,h)$, then $f$ is birational onto its image since $(\Gamma,h)$ is simple, and thus $[f(C)]\in V^\irr_{g,\Delta}$. 
    
    To finish the proof, it remains to show that $(\Gamma,h)$ admits a unique lift to a parametrized curve $f\:C\to S_\Delta$ passing through the $p_i$'s. Indeed, the uniqueness assertion of the theorem will then follow from the choice of the $p_i$'s as in Proposition~\ref{prop:lifting genral curve through points to all comp}, and the existence from Lemma~\ref{lem:surjectivity into nice cones} since $(\Gamma,h)\in M_\Theta$ is general. To do so, we will use Theorem~\ref{thm:corshort}.
    
    First, let us enhance the curve $(\Gamma,h)$ to a curve $(\Gamma',h')$ by adding a vertex $w_i$ and a contracted leg $l_i$ at each $h^{-1}(q_i)\in\Gamma$. As usual, we shift the labeling so that the contracted legs are the first $n$ legs. Since $(\Gamma,h)$ is simple, $(\Gamma',h')$ clearly satisfies condition (1) of Theorem~\ref{thm:corshort}, and condition (2) is satisfied by the construction. 
    
    To verify condition (3) and to conclude the uniqueness, it remains to show that $|\det(\zeta)|=1$. By \cite[Lemma~4.20]{Mik05}, the complement of $\{w_i\}_{i=1}^n$ is a union of trees and legs, and each tree $\Gamma_j$ contains a unique leg $l_j$ with $j>n$. Set $w_j:=\ft(l_j)$. Since changing the orientation of the edges preserves $|\det(\zeta)|$, we may assume that for any $j$, the restriction of the orientation on $\Gamma$ to $\Gamma_j$ coincides with the orientation of a rooted tree with root $w_j$. To prove that $|\det(\zeta)|=1$, it remains to show that the complex $L^\bullet(\Gamma')$ is quasi-isomorphic to the trivial one. 

    Pick one of the trees $\Gamma_j$, and consider the two edges $e,e'$ adjacent to $w_j$. Consider the (possibly reducible) tropical curve $\tilde{h}\:\tilde{\Gamma}\to N_\RR$ obtained from $\Gamma'$ by removing $w_j$ and $l_j$, and replacing $e$ and $e'$ with legs of the same slopes. Notice that the connected components of $(\tilde{\Gamma},\tilde{h})$ are again simple, passing through the $\{q_i\}_{i=1}^n$, and the complement of $\{w_i\}_{i=1}^n$ is a union of trees and half-lines, while each tree contains a unique leg. 
        
    We have a natural surjection $L^\bullet(\Gamma')\to L^\bullet(\tilde{\Gamma})$, whose kernel is 
    $$\ZZ^4=N_{w_j}\times \ZZ_e\times\ZZ_{e'}\to N_e\times N_{e'}=\ZZ^4,$$ 
    where the map is given by the matrix 
    $\begin{pmatrix}
       I_2 & \frac{\partial h}{\partial \vec e} & 0  \\
       I_2 & 0 & \frac{\partial h}{\partial \vec e'}  
    \end{pmatrix}.$
    The absolute value of the determinant of this matrix is given by $\left|\det \left(\frac{\partial h}{\partial \vec e},\frac{\partial h}{\partial \vec e'} \right)\right|,$ which is equal to 1 since $(\Gamma, h)$ is a multiplicity-one curve. Therefore, $L^\bullet(\Gamma')\to L^\bullet(\tilde{\Gamma})$ is a quasi-isomorphism. Proceeding this way until the only vertices left are the $w_i$'s, we obtain a quasi-isomorphism from $L^\bullet(\Gamma')$ to the complex $\prod_{i=1}^n N_{w_i}\to N^n$ with the identity map. Thus, $L^\bullet(\Gamma')$ is quasi-isomorphic to the trivial complex as asserted.
\qed

\begin{rem}
    It is well-known that for a simple tropical curve, Mikhalkin's multiplicity is given by the absolute value of the determinant of the evaluation map from the tropical moduli space to $N^n$; see, e.g., Hannah Markwig's thesis \cite[Lemma~4.68]{Mar06}. One can show that the two-term complex of the evaluation map is quasi-isomorphic to $L^\bullet(\Gamma')$, which explains why $|\det(\zeta)|$ should be 1. In fact, a straightforward modification of the argument we used shows directly that Mikhalkin's multiplicity of a simple parametrized tropical curve is given by $|\det(\zeta)|$.  
\end{rem}

\section{The combinatorics of simple stretched floor decomposed strata} \label{sec:algorithm}

In this section, we establish Theorem~C, the strong connectedness result for the tropical moduli space $M_{g, \nabla}^\trop$ for most tropical degrees $\nabla$ dual to $h$-transverse polygons. Throughout the section, $M=\ZZ^2$. 

\begin{defn}\label{def:swequiv}
    We say that nice strata $M_{[\Theta]}, M_{[\Theta']}\subset M_{g,\nabla}^\trop$ are {\em sw-equivalent} (or {\em connected thro\-ugh simple walls}) if there exists a sequence $M_{[\Theta]}=M_{[\Theta_1]},M_{[\Theta_2]},\dotsc,M_{[\Theta_k]}=M_{[\Theta']}$ of nice strata in $M_{g,\nabla}^\trop$ such that $M_{[\Theta_i]}$ and $M_{[\Theta_{i+1}]}$ share a simple wall as a common face. Again, by a slight abuse of language, we sometimes say that $M_\Theta$ and $M_{\Theta'}$ are sw-equivalent if $M_{[\Theta]}$ and $M_{[\Theta']}$ are.
\end{defn}

\begin{ThmC}
    Let $\Delta\subset \RR^2$ be an $h$-transverse polygon and $\nabla$ a tropical degree dual to $\Delta$. Suppose the non-vertical slopes in $\nabla$ are all primitive and $\nabla$ contains either $(0,1)$ or $(0,-1)$. Then any two ssfd strata in $M_{g,\nabla}^\trop$ are sw-equivalent. 
\end{ThmC}

The rest of the section is devoted to the proof of Theorem~C. In \S~\ref{subsec:basic moves}, we discuss two modifications of floor decomposed tropical curves: (1) moving an elevator along a floor, and (2) going up/down an elevator. Both modifications realize the sw-equivalence of the corresponding strata. In \S~\ref{subsec:self intersection}, we show that any ssfd stratum is sw-equivalent to an ssfd stratum $M_{[\Theta]}$ containing an ssfd curve without self-intersections. In \S~\ref{subsec:multiplicity sequence}, we introduce the multiplicity sequence for ssfd curves without self-intersections and describe explicitly the minimal multiplicity sequences. This is the step where we use the assumption about the existence of a leg-elevator of multiplicity one (i.e., that $\nabla$ contains either $(0,1)$ or $(0,-1)$). Finally, in \S~\ref{subsec:proofThmC}, we prove Theorem~C by using the results of \S~\ref{subsec:self intersection} and \S~\ref{subsec:multiplicity sequence} and showing that any pair of strata containing curves with minimal multiplicity sequences is sw-equivalent. We often slightly abuse terminology and switch between edges/legs and their images; for example, we say ``the $x$-coordinate of an elevator $E$'' and write $x(E)$ instead of ``the $x$-coordinate of $h(E)$''.

\subsection{The basic moves} \label{subsec:basic moves}

One of the advantages of floor decomposed curves is that their local deformations are particularly easy to describe: they are given by moving floors vertically and elevators horizontally; see, e.g., \cite[Proposition 5.9]{BM08}. Since this is central to our constructions below, we include a formal statement and a sketch of the proof for the reader's convenience:

\begin{lem}\label{lem:ssfddeform} 
    Let $(\Gamma_0, h_0)$ be a simple floor decomposed curve. Denote its floors by $F_1,\dotsc, F_\tth$, its elevators by $E_1,\dotsc, E_l$, and its combinatorial type by $\Theta_0$. Pick a vertical line $L$. For any curve $(\Gamma,h)$ in a neighborhood of $(\Gamma_0, h_0)$ in $M_{\Theta_0}$, and for each $1\le j\le \tth$, denote by $q_j$ the point of intersection of $L$ with $h(F_j)$. Set $y_j(\Gamma,h):=y(q_j)$ and $x_i(\Gamma,h):=x(h(E_i))$. Then $((x_i(\Gamma,h))_i,(y_j(\Gamma,h))_j)$ is a system of coordinates in a neighborhood of $(\Gamma_0, h_0)$ in $M_{\Theta_0}$.
\end{lem}

\begin{proof}
    First, assume that $\tth=1$. Since $\Gamma_0$ is 3-valent and has a unique floor, its elevators have pairwise distinct $x$-coordinates. After relabeling the elevators if needed, we may assume that $x_1(\Gamma_0,h_0)<\dots<x_l(\Gamma_0,h_0)$. We claim that the map $(\Gamma,h)\mapsto ((x_i(\Gamma,h))_i,y_1(\Gamma,h))$ identifies $M_{\Theta_0}$ with the cone $C_0$ given by $x_1<\dots<x_l$. Indeed, since $M_{\Theta_0}$ parametrizes simple curves, the image belongs to the complement of the $x$-diagonals, and since $M_{\Theta_0}$ is connected, the image belongs to $C_0$. Conversely, since all edges belong to the floor and have non-zero slopes, for any point $((x_i)_i,y_1)\in C_0$, the equations $x_i(\Gamma,h)=x_i$ for $i=1,\dotsc, l$ determine the lengths of the edges and the parametrization $h\:\Gamma\to\RR^2$ up to a vertical translation. Together with $y_1(\Gamma,h)=y_1$, they uniquely determine the curve $(\Gamma,h)\in M_{\Theta_0}$.

    The general case follows from the case $\tth=1$ because the curve $(\Gamma_0, h_0)$ is glued from $\tth$ curves, each having a single floor, along the elevators. More precisely, by cutting the bounded elevators and replacing each of them with a pair of elevator legs, we obtain $\tth$ simple floor decomposed curves $(\Gamma_j, h_j)$ of types $\Theta_j$. If we identify each $M_{\Theta_j}$ with the corresponding cone $C_j$, then a small neighborhood of $(\Gamma_0, h_0)\in M_{\Theta_0}$ can be identified with the fibered product of neighborhoods of $(\Gamma_j, h_j)\in C_j$ along suitable projections. For a pair $j_1,j_2$, the projections are given by the $x$-coordinates of the elevator legs corresponding to the common elevators adjacent to the floors $F_{j_1}$ and $F_{j_2}$ in the curve $(\Gamma_0, h_0)$. Since $((x_i)_i,(y_j)_j)$ is a system of coordinates on the fibered product, the lemma follows.
\end{proof}

\begin{rem}
    If the curve $(\Gamma_0,h_0)$ is stretched, then, similarly to the case $\tth=1$, the coordinates $((x_i(\Gamma,h))_i,(y_j(\Gamma,h))_j)$ are well-defined beyond a small neighborhood of $(\Gamma_0,h_0)$. For example, if $a<x_1(\Gamma_0,h_0)<\dots<x_l(\Gamma_0,h_0)<b$ and the difference of the $y$-coordinates of any pair of vertices belonging to different floors is much larger than $b-a$, then the coordinate chart contains the simplex given by $a<x_1<\dots<x_l<b$ and $y_j=y_j(\Gamma_0,h_0)$ for all $j$.
\end{rem}
    
In the following constructions, we work with an ssfd curve $(\Gamma_0, h_0)$ whose elevators have pairwise distinct $x$-coordinates. Note that by Lemma~\ref{lem:ssfddeform}, any ssfd curve admits a small perturbation satisfying this property. Consider the ``left to right'' labeling of the elevators, i.e., $i>j$ if and only if $x(h_0(E_i))>x(h_0(E_j))$; similarly, consider the ``bottom to top'' labeling of the floors $F_1, \dotsc, F_{\tth}$. For convenience, we add {\em virtual floors} at infinity, $F_0$ and $F_{\tth +1}$, at heights $y=-\infty$ and $y=\infty$, respectively. We virtually connect these floors to the downward and upward leg-elevators, respectively. Unless excluded explicitly, we include the virtual floors whenever we address floors. 

\begin{figure}[ht]   
    \begin{tikzpicture}[x=0.8pt,y=0.8pt,yscale=-.9,xscale=.9]
        \import{./}{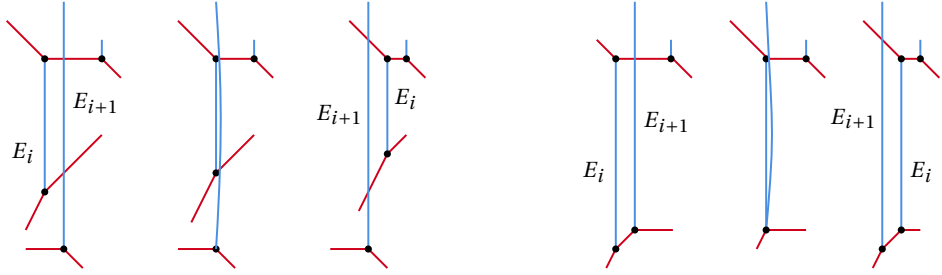}
    \end{tikzpicture}   
    \caption{Two illustrations of Construction~\ref{constr:moving right}: on the left, the elevators $E_i$ and $E_{i+1}$ are not adjacent to the same non-virtual floor, and on the right, they are. In both pictures, the curve $(\Gamma_0,h_0)$ is on the left, and $(\Gamma_1,h_1)$ is in the middle.}
    \label{fig:moving right}
\end{figure}

\begin{constr}[Changing the relative order of elevators]\label{constr:moving right}
    Assume that the elevators $E_i$ and $E_{i+1}$ are either both unbounded or do not connect the same pair of floors. Pick a vertical line $L$, and in the notation of Lemma~\ref{lem:ssfddeform}, consider the deformation of $(\Gamma_0,h_0)$ given by fixing all $y$-coordinates and all $x$-coordinates except $x_i$, and increasing the coordinate $x_i$. Since $(\Gamma_0,h_0)$ is stretched, increasing $x_i$ eventually leads to a stretched floor decomposed curve $(\Gamma_1,h_1)$ in the closure $\overline{M}_{\Theta_0}$ for which $x_i=x_{i+1}$. If $E_i$ and $E_{i+1}$ are not adjacent to the same non-virtual floor, then $(\Gamma_1,h_1)\in M_{\Theta_0}$, and increasing $x_i$ further switches the relative order of $E_i$ and $E_{i+1}$. Otherwise, the curve $(\Gamma_1,h_1)$ belongs to a simple wall $M_{\Theta_1}$ by our assumption on the elevators $E_i$ and $E_{i+1}$. One of the other two nice strata adjacent to $M_{\Theta_1}$ then contains an ssfd deformation of $(\Gamma_1,h_1)$ for which $x_i>x_{i+1}$, while all coordinates except $x_i$ remain fixed; see Figure~\ref{fig:moving right} for an illustration. By construction, this stratum is sw-equivalent to $M_{\Theta_0}$. 
\end{constr}

\begin{constr}[Going up/down an elevator]\label{constr:moving up}
    Assume that the elevator $E_i$ is adjacent to the floors $F_k$ and $F_s$, and the elevator $E_{i+1}$ to $F_k$ and $F_r$. Suppose that one of the following inequalities holds: $1 \leq k < s < r\le \tth+1$, or $1 \leq k < r < s\le \tth+1$, or $0 \leq s < r < k\le \tth$, or $0 \leq r < s < k\le \tth$. The construction in all four cases is essentially the same; thus, we describe it only in the first case.
    
    Again, let $L$ be a vertical line. Since the curve is stretched and $k < s < r$, it follows that $h_0(E_{i+1})$ intersects $h_0(F_s)$. Denote the preimage of the intersection point on $F_s$ by $p$, and consider the tropical curve $\Gamma'$ obtained from $\Gamma_0$ by disconnecting the elevator $E_{i+1}$ from the floor $F_k$ and attaching it to the floor $F_s$ at the point $p$ instead. Then there exists an ssfd stratum $M_{\Theta'}$ sw-equivalent to $M_{\Theta_0}$ and a parametrization $h'\:\Gamma'\to\RR^2$ such that $(\Gamma',h')$ is an ssfd curve whose intersection with $L$ coincides with that of $\Gamma_0$ floorwise. Furthermore, the $x$-coordinates of the elevators of $\Gamma'$ coincide with those of the elevators of $\Gamma_0$, and the slopes (including multiplicities) along all common edges and legs\footnote{Note that there is an edge in $\Gamma_0$ that splits into two, and a pair of edges are glued into one in $\Gamma'$.} between $h_0$ and $h'$ coincide, except for $E_i$. Moreover, $\fw(h',E'_i)=\fw(h_0, E_i)+\fw(h_0, E_{i+1})$, where we denote the elevator $E_i$ in $\Gamma'$ by $E'_i$; cf. Figure~\ref{fig:moving up}.

    \begin{figure}[ht]
    \begin{tikzpicture}[x=0.8pt,y=0.8pt,yscale=-.9,xscale=.9]
        \import{./}{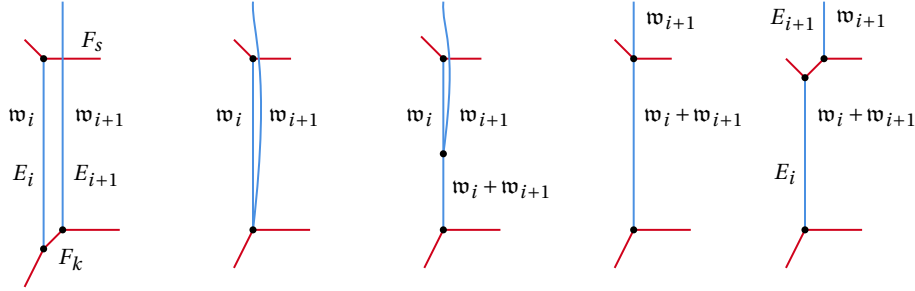}
    \end{tikzpicture}   
    \caption{An illustration of Construction~\ref{constr:moving up}. On the left is the curve $(\Gamma_0,h_0)$; on the right is $(\Gamma',h')$. The curves in the intermediate strata $M_{\Theta_1}, M_{\Theta_2}, M_{\Theta_3}$ are depicted in the middle. The multiplicities at each step are denoted by $\fw_\bullet$.}
    \label{fig:moving up}
    \end{figure}
    
    To see this, we again use the coordinates from Lemma~\ref{lem:ssfddeform}. Consider the deformation of $(\Gamma_0,h_0)$ induced by fixing the $y_j$'s and all the $x_j$'s except $x_i$, while increasing $x_i$. As in Construction~\ref{constr:moving right}, there exists a curve $(\Gamma_1,h_1)\in \overline{M}_{\Theta_0}$ for which $x_i=x_{i+1}$. As before, $(\Gamma_1,h_1)$ belongs to a simple wall $M_{\Theta_1}$ because $r\ne s$. Since $r,s>k$, there exists a nice (but not simple) stratum $M_{\Theta_2}$ adjacent to $M_{\Theta_1}$ corresponding to the resolution of the 4-valent vertex, in which the two elevators remain adjacent to the same vertex $u$; see Figure~\ref{fig:moving up} for an illustration. By the balancing condition, for any curve $(\Gamma_2,h_2)\in M_{\Theta_2}$, the outgoing slopes at $u$ sum to zero; hence, $u$ must lie above $F_k$. Moreover, the original elevator $E_i$ in $M_{\Theta_1}$ is split by $u$ into two elevators in $\Gamma_2$. We denote the one adjacent to $u$ from below by $E'_i$, and the one from above by $E_i$. Again, by the balancing condition at $u$, we have 
    $$ \fw(h_2, E'_i)=\fw(h_2, E_i)+\fw(h_2, E_{i+1})=\fw(h_0, E_i)+\fw(h_0, E_{i+1}). $$ 
    By moving the vertex $u$ up and adjusting the lengths of the elevators $E_i, E_{i+1}, E'_i$ accordingly, we deform the curve $(\Gamma_2,h_2)$ to a curve $(\Gamma_3,h_3)$, which corresponds to shrinking the length of $E_i$ to 0. The curve $(\Gamma_3,h_3)$ belongs to another simple wall $M_{\Theta_3}$ and has a 4-valent vertex on floor $F_s$, which we denote by $v$. By construction, there are two elevators adjacent to $v$: $E'_i$ and $E_{i+1}$. Finally, there is a nice stratum $M_{\Theta'}$ adjacent to $M_{\Theta_3}$ in which the 4-valent vertex $v$ is resolved in such a way that $x_i<x_{i+1}$; see Figure~\ref{fig:moving up} again. The stratum $M_{\Theta'}$ contains the desired curve $(\Gamma',h')$.
\end{constr}

\begin{figure}[ht]
    \begin{tikzpicture}[x=0.8pt,y=0.8pt,yscale=-.9,xscale=.9]
        \import{./}{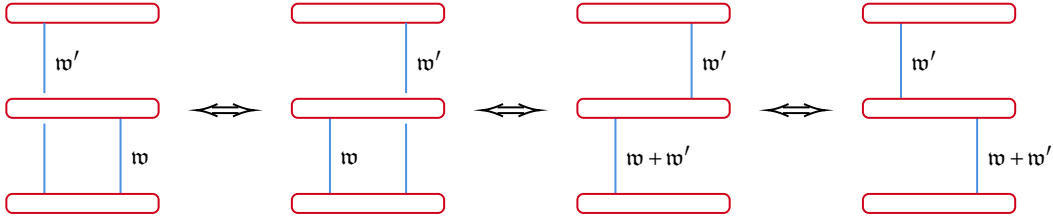}
    \end{tikzpicture}   
    \caption{A diagram visualizing the change of order of elevators and going up/down an elevator in the case $k<r,s$. In the case $k>r,s$, the diagram is simply the reflection of the diagram above with respect to the $x$-axis.}
    \label{fig:diagrams}
\end{figure}

The schematic visualization in Figure~\ref{fig:diagrams} of the equivalences induced by the two constructions will be useful in the proof of Theorem~C.

\subsection{Resolving self-intersections}\label{subsec:self intersection}

\begin{defn}
    Let $ M_{[\Theta]}$ be an ssfd stratum. The {\em self-intersection number} of $ M_{[\Theta]}$ is the minimal number of self-intersection points of an ssfd curve $(\Gamma, h)\in M_{\Theta}$.
\end{defn}

\begin{rem}
    For different ssfd curves in the same stratum, the number of self-intersection points may vary. For example, consider the left and middle curves in Figure~\ref{fig:comtype}. Although the middle curve is not ssfd, there is an ssfd curve in its stratum with the same combinatorics of the image.
\end{rem}

\begin{lem} \label{lem:self intersection}
    Any ssfd stratum $M_{[\Theta]}$ is sw-equivalent to an ssfd stratum with self-intersection number 0.
\end{lem}

\begin{proof}
    Let $(\Gamma,h)\in M_{\Theta}$ be an ssfd curve with at least one self-intersection point. To prove the lemma, it suffices to construct an sw-equivalent stratum $M_{\Theta'}$ and an ssfd curve $(\Gamma',h')\in M_{\Theta'}$ having a smaller number of self-intersection points. Because the curve is stretched, there are only two possible types of self-intersection points. The first is an intersection between an elevator $E$ and a floor $F$, which occurs whenever $E$ connects a floor located above $F$ to one located below it. The second is an intersection between two floors, which can only occur between non-elevator legs whose slopes share the same $x$-component (either both 1 or both $-1$). More precisely, if $k<s$ and the legs $l_k\subset F_k$, $l_s\subset F_s$ have slopes $(\pm 1,y_k)$, $(\pm 1,y_s)$ sharing the same $x$-component, then $h(l_k)\cap h(l_s)\ne\emptyset$ if and only if $y_k>y_s$. We will prove the existence of the desired curve $(\Gamma',h')$ by considering these two types separately and applying the constructions from \S~\ref{subsec:basic moves}. 
    
    After a small perturbation, we may assume that the $x$-coordinates of the elevators of $(\Gamma,h)$ are pairwise distinct. As usual, we label the elevators ``left to right,'' the floors ``bottom to top,'' and add the virtual floors $F_0$ and $F_{\tth +1}$. Below, we will say that two floors $F$ and $F'$ are {\em adjacent} if they are joined by an elevator.
    
    Case 1: \emph{There exists a self-intersection point between a floor and an elevator.} We claim that there exists $1\le k\le \tth$ such that the floor $F_k$ is adjacent to two distinct floors, either both above $F_k$ or both below $F_k$. Assuming the claim, we can construct the desired curve $(\Gamma',h')$ as follows. Assume without loss of generality that the floor $F_k$ is adjacent to two distinct floors above it; i.e., there exist $k<s<r$ such that $F_k$ is adjacent to $F_s$ and $F_r$ (the argument in the second case is similar). Let $E_i$ and $E_j$ be elevators connecting $F_k$ to $F_s$ and to $F_r$, respectively. After applying Construction~\ref{constr:moving right} several times, we may assume that $j=i+1$. Note that applying this construction to an ssfd curve preserves the number of self-intersection points of either type. Now, we can apply Construction~\ref{constr:moving up}. The resulting ssfd curve $(\Gamma',h')$ has fewer self-intersection points than $(\Gamma,h)$ because the intersection point between $E_{i+1}$ and $F_s$ is resolved; see Figure~\ref{fig:moving up} for an illustration. By construction, the stratum of $(\Gamma',h')$ is sw-equivalent to $M_{\Theta}$, and we are done.

    Let us now prove the existence of $k$. Suppose an elevator $E$ intersects a floor $F$. It must connect a floor above $F$ to a floor below it; let us denote these floors by $F_s$ and $F_{s'}$, respectively. We will assume that $F_s$ is not virtual (at least one of the floors $F_s$ and $F_{s'}$ is not, and the argument in the other case is similar). Because $\Gamma$ is connected, there is a shortest sequence of non-virtual floors $F_{s}=F_{r_1},F_{r_2},\dotsc,F_{r_d}=F_{s-1}$ such that $F_{r_i}$ is adjacent to $F_{r_{i+1}}$ for any $1\leq i\leq d-1$. If $d = 2$, then $k=s$ satisfies the assertion, as $F_s$ is adjacent to $F_{s-1}$ and to $F_{s'}$. Otherwise, either $\max\{r_i\}>s$ or $\min\{r_i\}<s-1$. In the former case, $k=\max\{r_i\}$ satisfies the assertion, and in the latter, $k=\min\{r_i\}$ does.

    Case 2: \emph{All self-intersection points are between pairs of floors.} Let $q$ be a self-intersection point between the floors $F_k$ and $F_s$, where $1\le k<s\le\tth$. Since $(\Gamma, h)$ is stretched, $q$ belongs to a pair of non-elevator legs: $l_k$ of $F_k$ and $l_s$ of $F_s$. Furthermore, the slopes of these legs share the same $x$-component. We will assume that this component is equal to $-1$, as the argument in the second case is similar. We claim that we can choose $q$ such that $s=k+1$. Indeed, assume $s>k+1$. Since the curve is stretched, $h$ maps the vertices of $F_{k+1}$ strictly above $h(F_k)$ and below $h(F_s)$. Let $l_{k+1}$ be the leg of $F_{k+1}$ whose slope has $x$-component $-1$. It follows that $h(l_{k+1})$ must intersect $h(l_k)\cup h(l_s)$. If $h(l_{k+1})$ intersects $h(l_k)$, then the intersection point $h(l_k)\cap h(l_{k+1})$ serves as the desired replacement for $q$. Otherwise, we replace $k$ with $k+1$, replace $q$ with $h(l_{k+1})\cap h(l_s)$, and proceed inductively.

    Let us now construct the desired curve $(\Gamma',h')$. Since $(\Gamma,h)$ is connected and has no self-intersection points between elevators and floors, any two consecutive floors must be adjacent to each other. Let $E_i$ be the elevator with the minimal index connecting $F_k$ and $F_{k+1}$. Choose $a\gg 1$ and let $L$ be the line given by $x=a$. Consider the coordinates on $M_\Theta$ as in Lemma~\ref{lem:ssfddeform}, and the deformation of $(\Gamma,h)$ given by fixing all coordinates except $x_i$, and decreasing $x_i$. Eventually, this leads to a curve $(\Gamma_1,h_1)$ belonging to a simple wall $M_{\Theta_1}$ in which the length of $E_i$ shrinks to zero. 

    \begin{figure}[ht]   
    \begin{tikzpicture}[x=0.8pt,y=0.8pt,yscale=-.9,xscale=.9]
        \import{./}{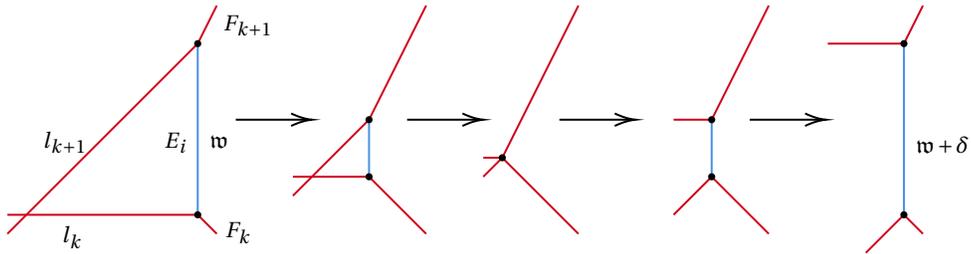}
    \end{tikzpicture}   
    \caption{An illustration of the reduction of the self-intersection number in the proof of Lemma~\ref{lem:self intersection}.}
    \label{fig:wc3}
    \end{figure}
    
    The simple wall $M_{\Theta_1}$ is adjacent to another ssfd stratum $M_{\Theta'}$ in which the leg $l_k$ is attached to the floor $F_{k+1}$ and $l_{k+1}$ to $F_k$; see Figure~\ref{fig:wc3} for an illustration. Since the $x$-components of the slopes of $l_k$ and $l_{k+1}$ are equal to $-1$, the new edge $E'_i$ in $\Theta'$ is an elevator, and the resolution is again simple and floor decomposed. By the balancing condition, $\fw(h',E'_i)=\fw(h,E_i)+\delta$, where $\delta$ is the absolute value of the difference between the $y$-coordinates of the slopes of $l_k$ and $l_{k+1}$. Let us deform the curve $(\Gamma_1,h_1)$ into the stratum $M_{\Theta'}$ by fixing all coordinates except $x_i$ and increasing $x_i$ until $x_i(\Gamma',h')=x_i(\Gamma,h)$. The resulting curve $(\Gamma',h')$ is an ssfd curve. Its number of self-intersections is one fewer than that of $(\Gamma,h)$ because the point $q$ is resolved, while other intersection points between floors are preserved and no elevator-floor intersection points are produced. Finally, by construction, $M_{\Theta'}$ is sw-equivalent to $M_{\Theta}$, which completes the proof.
\end{proof}

\subsection{The multiplicity sequence}\label{subsec:multiplicity sequence} 

Let $(\Gamma_0,h_0)$ be an ssfd curve of degree $\nabla$ and genus $g$ without self-intersections, and let $M_{\Theta_0}$ be its stratum. Since $\Gamma_0$ is connected, the floor $F_r$ is connected to the floors $F_{r-1},F_{r+1}$ and no other floors for any $1<r<\tth$. Furthermore, the virtual floor $F_0$ (resp., $F_{\tth +1}$) may only be connected to $F_1$ (resp., $F_\tth$). Denote the upward elevators adjacent to the floor $F_r$ with $0\le r\le \tth$ by $E_{r,i}$ for $1\le i\le j_r$, where the indices are chosen such that $i_1<i_2$ if and only if $x(E_{r,i_1})<x(E_{r,i_2})$, and set $\fw_{r,i}:=\fw (E_{r,i})$.

\begin{defn}\label{def:multiplicity sequence}
    The {\em multiplicity sequence} of $(\Gamma_0,h_0)$ is the sequence $\left(\fw_{r,i}\right)_{0\le r\le\tth}^{1\le i\le j_r}$ ordered lexicographically with respect to the pairs of indices, i.e., $\fw_{0,1},\fw_{0,2},\dotsc, \fw_{\tth,j_\tth}$. Denote the length of the multiplicity sequence by $\ell(\Gamma_0,h_0):=\sum_{r=0}^\tth j_r$.
\end{defn}

\begin{rem}\label{rem:multseqprop}
    (1) The length $\ell(\Gamma_0,h_0)$ of the multiplicity sequence is uniquely determined by $\nabla$ and $g$. Indeed, the multiset $\{(0, -\fw_{0,i})\}_{1\le i\le j_0}$ is exactly the multiset of vertical vectors in $\nabla$ with a negative $y$-coordinate. Similarly, the multiset of vertical vectors with a positive $y$-coordinate is $\{(0, \fw_{\tth,i})\}_{1\le i\le j_\tth}$. Thus, $j_0$ and $j_\tth$ are determined by $\nabla$. Now, since $(\Gamma_0,h_0)$ is an ssfd curve, its genus is given by $g=\sum_{r=1}^{\tth -1}(j_r-1)$, which implies our claim.

    (2) Let $\Delta$ be an $h$-transverse polygon dual to $\nabla$. Denote by $a_r$ the length of the slice of $\Delta$ at height $y=r_0+r$, where $r_0:=\min_{(x,y)\in \Delta}y$. Then the multiplicity sequence necessarily satisfies the following system of equations: $\sum_{i=1}^{j_r}\fw_{r,i}=a_r$ for all $0\leq r\leq\tth$. In particular, the set of multiplicity sequences of degree-$\nabla$ ssfd curves without self-intersections is finite.
\end{rem}

\begin{lem}\label{lem:minmultseq}
    Let $(\Gamma_0,h_0)$ be an ssfd curve of degree $\nabla$ and genus $g$ without self-intersections, and let $M_{\Theta_0}$ be its stratum. Assume that $(0,-1)\in\nabla$ and that the multiplicity sequence of $(\Gamma_0,h_0)$ is minimal with respect to the lexicographic order among all multiplicity sequences of ssfd curves without self-intersections belonging to strata in the sw-equivalence class of $M_{\Theta_0}$. Then, there exists a maximal index $1\leq k\leq\tth-1$, uniquely determined by $\nabla$ and $g$, such that
    \begin{enumerate}
        \item $\fw(E_{r,i})=1$ for all $1\le r<k$ and $1\le i\le j_r$,
        \item $\fw(E_{k,i})=1$ for all $1\le i<j_k$, and
        \item $j_r=1$ and $\fw(E_{r,1})=a_r$ for all $k<r<\tth$.
    \end{enumerate}
\end{lem}

\begin{proof}
    We start with the uniqueness of $k$. By Remark~\ref{rem:multseqprop} (1), the genus $g$ of the ssfd curve $\Gamma_0$ is given by the number of bounded elevators minus the number of floors plus 1. Thus, the integer $k$ satisfying the conclusion of the lemma is either the minimal integer for which $\sum_{r=1}^k (a_r-1)>g$, or $\tth-1$ if no such integer exists. Since the $a_r$'s are uniquely determined by $\Delta$, which is determined by $\nabla$, the integer $k$ is uniquely determined by $\nabla$ and $g$. Let us now prove the existence of $k$.
    
    If $\fw(E_{r,i})=1$ for all $1\le r\le \tth-1$ and all $1\le i\le j_r$, then $k:=\tth-1$ clearly satisfies the assertion. Otherwise, let $1\leq k\leq\tth-1$ be the minimal index for which there exists an $i$ such that $\fw(E_{k,i})\ge 2$, and let $j$ be the minimal such $i$. It is sufficient to show that $j=j_k$ and $j_r=1$ for all $k<r<\tth$. 
    
    Step 1: Let us show that $j=j_k$. Assume to the contrary that $j<j_k$. Let $E_{k-1,i}$ be an elevator of multiplicity one. Such an elevator exists due to our assumption on the degree $\nabla$ and the minimality of $k$. Since Construction~\ref{constr:moving right} does not change the multiplicity sequence or the sw-equivalence class, after applying it several times, we may assume that 
    \[ x(E_{k-1,i-1})<x(E_{k,j})<x(E_{k-1,i})<x(E_{k,j+1})<x(E_{k-1,i+1}). \] 
    Now, applying the inverse of Construction~\ref{constr:moving up} to the pair $E_{k,j}, E_{k-1,i}$, and then Construction~\ref{constr:moving up} to the pair $E_{k-1,i}, E_{k,j+1}$, we obtain an ssfd curve with a smaller multiplicity sequence than that of $(\Gamma_0,h_0)$, which is a contradiction; see Figure~\ref{fig:mult1}. Thus, $j=j_k$. 
    
    \begin{figure}[ht]   
        \begin{tikzpicture}[x=0.8pt,y=0.8pt,yscale=-0.9,xscale=0.9]
            \import{./}{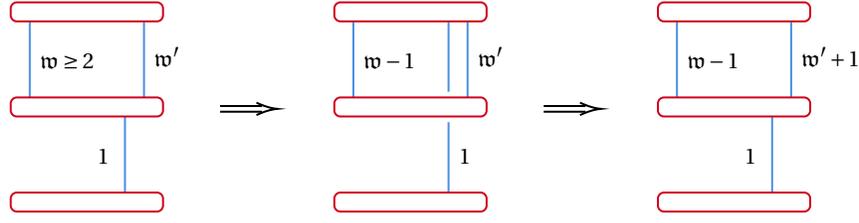}
        \end{tikzpicture}   
        \caption{The application of Construction~\ref{constr:moving up} to reduce the multiplicity sequence in Step~1. Here the diagrams represent the floors $F_{k-1},F_k,F_{k+1}$, the elevators $E_{k-1,i}, E_{k,j},E_{k,j+1}$, and the multiplicities $\fw:=\fw(E_{k,j})$, $\fw':=\fw(E_{k,j+1})$, and $1=\fw(E_{k-1,i})$.}
        \label{fig:mult1}
    \end{figure}

    Step 2: Let us show that $j_r=1$ for all $k<r<\tth$. If not, pick the minimal $k<r<\tth$ such that $j_r>1$. Assume first that there exist $E_{r,j}$ and $E_{r-1,i}$ such that $\fw(E_{r-1,i})>\fw(E_{r,j})$. Without loss of generality, we may assume that $j<j_r$. As in Step~1, we may assume that 
    \[ x(E_{r-1,i-1})<x(E_{r,j})<x(E_{r-1,i})<x(E_{r,j+1})<x(E_{r-1,i+1}). \]
    Now, applying the two constructions and their inverses as in Figure~\ref{fig:mult2}, we obtain an ssfd curve with a smaller multiplicity sequence than that of $(\Gamma_0,h_0)$, whose stratum is sw-equivalent to $M_{\Theta_0}$, which is a contradiction. Therefore, 
    \begin{equation}\label{eq:multseqredin}
        \fw(E_{r-1,i})\le\fw(E_{r,j})
    \end{equation}
    for all $1\le j\le j_r$ and $1\le i\le j_{r-1}$.

    \begin{figure}[ht]   
        \begin{tikzpicture}[x=0.8pt,y=0.8pt,yscale=-0.9,xscale=0.9]
            \import{./}{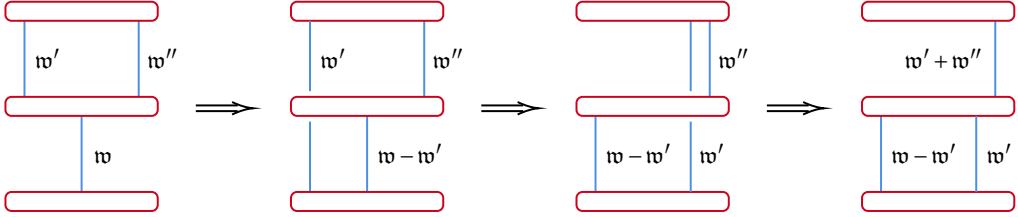}
        \end{tikzpicture}   
        \caption{The first reduction of the multiplicity sequence in Step~2. Here the diagrams represent the floors $F_{r-1},F_r,F_{r+1}$, the elevators $E_{r-1,i}, E_{r,j},E_{r,j+1}$, and the multiplicities $\fw:=\fw(E_{r-1,i})$, $\fw':=\fw(E_{r,j})$, and $\fw'':=\fw(E_{r,j+1})$.}
        \label{fig:mult2}
    \end{figure}
    
    We claim that $r=k+1$ and $\fw(E_{k,1})=1$. To see this, it is sufficient to show that $j_{r-1}>1$. Indeed, if either $r>k+1$ or $\fw(E_{k,1})\ne 1$, then $j_{r-1}=1$ by the choice of $k$ and the minimality of $r$. Assume that $j_{r-1}=1$. Then $a_{r-1}=\fw(E_{r-1,1})\le \frac{1}{2}a_r$ by \eqref{eq:multseqredin} since $j_r>1$. However, $2a_{r-1}\ge a_r+a_{r-2}$ by the convexity of $\Delta$. Therefore, $a_{r-2}=0$, which is a contradiction since $r>k\ge 1$ and $\Delta$ has a flat bottom side by the assumption on $\nabla$. Note that $\fw(E_{r,j})\ge\fw(E_{k,j_k})\ge 2$ for all $1\le j\le j_r$ by \eqref{eq:multseqredin}. Applying the constructions exactly as in Step~1 to the elevators $E_{r,1}, E_{r,2}$ and $E_{r-1,1}=E_{k,1}$, we again obtain an ssfd curve with a smaller multiplicity sequence than that of $(\Gamma_0,h_0)$, whose stratum is sw-equivalent to $M_{\Theta_0}$, which is a contradiction. The proof is now complete.
\end{proof}

\subsection{Proof of Theorem~C}\label{subsec:proofThmC}

Theorem~C follows rather easily from what we have proved so far. To apply Lemma~\ref{lem:minmultseq}, we will assume that $(0,-1)\in\nabla$, which can be achieved by reflecting $\Delta$ with respect to the $x$-axis if necessary. 

Let $M_{[\Theta]}\subset M_{g,\nabla}^\trop$ be any ssfd stratum and let $\Upsilon$ be its sw-equivalence class. By Lemma~\ref{lem:self intersection}, $\Upsilon$ contains an ssfd stratum with self-intersection number 0. Furthermore, we can choose an ssfd stratum $M_{[\Theta_0]}\in\Upsilon$ and an ssfd curve $(\Gamma_0,h_0)\in M_{\Theta_0}$ without self-intersections, whose multiplicity sequence is minimal among all ssfd curves without self-intersections in $\bigcup_{M_{[\Xi]}\in\Upsilon}M_{[\Xi]}$. The existence of such a minimizing curve is clear because the set of multiplicity sequences of curves of genus $g$ and degree $\nabla$ is finite; see Remark~\ref{rem:multseqprop}. Since Construction~\ref{constr:moving right} does not change the multiplicity sequence, after applying it several times, we may assume that the elevators of $(\Gamma_0,h_0)$ satisfy the following: $x(E_{r,i})\le x(E_{r+1,j})$ for all $0\le r\le \tth$ and all $i,j$.

It remains to prove that the stratum $M_{[\Theta_0]}$ is independent of $\Theta$ and depends only on $\nabla$ and $g$. To show this, note that the combinatorial type $\Theta_0$ is uniquely determined by the following data: (a) the slopes of the pairs of legs adjacent to each floor, (b) the sequences of multiplicities of the upward elevators adjacent to each floor, and (c) the relative positions of the $x$-coordinates of the upward elevators adjacent to the floors $F_r$ and $F_{r+1}$ for all $0\le r\le\tth$. 

First, since the ssfd curve $\Gamma_0$ has no self-intersections, the slopes of the legs at floor $F_r$ are dual to the non-horizontal sides of the slice $\Delta\cap \left(\RR\times[r_0+r, r_0+r+1]\right)$ and therefore depend only on $\Delta$, which is determined by $\nabla$. Second, since the number $k$ in Lemma~\ref{lem:minmultseq} depends only on $\nabla$ and $g$, the multiplicity subsequences in (b) do as well. Finally, by construction, $x(E_{r,i})\le x(E_{r+1,j})$ for all $0\le r\le \tth$ and all $i,j$. Therefore, the combinatorial type $\Theta_0$ and the stratum $M_{[\Theta_0]}$ depend only on $\nabla$ and $g$, which completes the proof. \qed

\section{Irreducibility of Severi varieties}\label{sec:Severi}

In this section, we establish the irreducibility of Severi varieties for a rich class of polarized toric surfaces (Theorem~B), namely for polarized toric surfaces associated to admissible polygons. We begin this section by defining admissible polygons and then proving Theorem~B. In \S~\ref{subsec:examples}, we discuss examples of such polygons. We provide an easily verifiable sufficient condition for an $h$-transverse polygon to be admissible (Proposition~\ref{prop:example of admissible}). In particular, we show the admissibility of the polygons of all classical polarized toric surfaces such as $\PP^2$, Hirzebruch surfaces, and toric del Pezzo surfaces equipped with arbitrary polarizations. Finally, we prove a version of Zariski's Theorem (Theorem~\ref{thm:zariski}) in \S~\ref{subsec:zarANDadj}.

\begin{defn} \label{def:admissible}
    A lattice polygon $\Delta \subset \RR^2$ is called {\em $g$-admissible} if 
    \begin{enumerate}
        \item $\Delta$ is $h$-transverse,
        \item $\Delta$ has a horizontal side, and
        \item there is an ssfd curve with multiplicity one in $M_{g,\nabla}^\trop$, where $\nabla$ is the reduced tropical degree dual to $\Delta$.
    \end{enumerate}
    We say that $\Delta$ is {\em admissible} if it is $g$-admissible for all $0 \leq g \leq |\Delta^\circ \cap \ZZ^2|$. 
\end{defn} 

As mentioned above, we give examples of admissible polygons in \S~\ref{subsec:examples}. Meanwhile, let us prove Theorem B, which follows from the following slightly more general statement:

\begin{thm}\label{thm:stongB}
    Suppose $\Delta$ is a $g$-admissible polygon. Then the Severi variety $V_{g, \Delta}^{\irr}$ is non-empty and irreducible over any algebraically closed field $k$.
\end{thm}

\begin{proof}[Proof of Theorem~\ref{thm:stongB}]
    To apply the tropical techniques developed in the previous sections, we must first reduce the assertion of the theorem to the case of a valued field $K$ as in \S~\ref{sec:not}. This is quite standard. Let $F:=k(\!(\varpi)\!)$ be the field of Laurent power series and $K$ be its algebraic closure. By \cite[Lemma 2.6]{CHT23}, the Severi varieties $V_{g,\Delta}^\irr(k)$ are locally closed subsets of the linear system $|\CL_\Delta|(k)$. The universal family of curves over $|\CL_\Delta|(K)$ is the base change of the universal family over $|\CL_\Delta|(k)$. Furthermore, the geometric genus of a curve is preserved under this base change. It follows that the locus $V_{g,\Delta}^\irr(K)\subset |\CL_\Delta|(K)$ is the base change of the one over $k$. Thus, from now on, we may assume that we work over $K$.

    Let $(\Gamma_0,h_0)$ be an ssfd multiplicity-one tropical curve of genus $g$ and reduced degree $\nabla$ dual to $\Delta$. Such a curve exists since $\Delta$ is admissible. Let $M_{[\Theta_0]}\subset M_{g,\nabla}^\trop$ be the stratum containing $(\Gamma_0,h_0)$. Let $V\subseteq V^\irr_{g,\Delta}$ be an irreducible component, and consider its tropicalization $\Sigma_V$ as defined in \S~\ref{subsec:tropseveri}. By Theorem~\ref{thm:existence and uniqueness}, there is a unique irreducible component $V_0\subseteq V^\irr_{g,\Delta}$ such that $M_{[\Theta_0]}\subset\overline{\Sigma}_{V_0}$. Thus, to prove the irreducibility of $V_{g, \Delta}^{\irr}$, it remains to show that $M_{[\Theta_0]}\subset\overline{\Sigma}_V$.
    
    Let $q_1,\dotsc, q_n\in\RR^2$ be vertically stretched points in general position, where $n:=\dim(V_{g,\Delta}^\irr)$, and let $p_1,\dotsc, p_n\in T\subset S_\Delta$ be algebraic points tropicalizing to the $q_i$'s. Then the $p_i$'s are in general position, and since the Severi variety is equidimensional, the locus of curves in $V$ passing through the $p_i$'s is zero-dimensional and non-empty. Let $[C]\in V$ be such a curve and $(\Gamma_1,h_1)$ its tropicalization. Then $h_1(\Gamma_1)$ contains the points $\{q_i\}$, and therefore $(\Gamma_1,h_1)$ is an ssfd curve by \cite[\S 5.1]{BM08}. Furthermore, since the $q_i$'s are general, it follows that $\dim(\Sigma_V \cap M_{[\Theta_1]})=\dim M_{[\Theta_1]}$, where $M_{\Theta_1}$ is the stratum containing $(\Gamma_1,h_1)$. Thus, $M_{[\Theta_1]}\subseteq \overline{\Sigma}_V$ by Lemma~\ref{lem:surjectivity into nice cones} (2). Finally, since $M_{[\Theta_1]}$ and $M_{[\Theta_0]}$ are sw-equivalent by Theorem~C, Lemma~\ref{lem:surjectivity into nice cones} (3) implies that $M_{[\Theta_0]}\subseteq \overline{\Sigma}_V$.
\end{proof}

\subsection{Examples of $g$-admissible polygons}\label{subsec:examples}

It is easier to work with a slightly more restrictive class of polygons, which we call {\em very admissible}. 

\begin{defn}\label{def:veryadmis}
    A {\em $g$-very admissible polygon} is an $h$-transverse polygon $\Delta\subset \RR^2$ with a horizontal side such that the moduli space $M_{g,\nabla}^\trop$ contains a unimodular ssfd curve, where $\nabla$ is the reduced tropical degree dual to $\Delta$. We say that the polygon $\Delta$ is {\em very admissible} if it is $g$-very admissible for all $0 \leq g \leq |\Delta^\circ \cap \ZZ^2|$. 
\end{defn} 
Since unimodular curves are multiplicity-one, any $g$-very admissible polygon is $g$-admissible. The following proposition provides a sufficient condition for a polygon $\Delta$ to be very admissible.

\begin{prop}\label{prop:example of admissible} 
    Let $\Delta$ be an $h$-transverse polygon with a horizontal side, and $a_r$ the length of the slice of $\Delta$ at height $y=r+r_0$, where $r_0:=\min_{(x,y)\in \Delta}y$. Assume that $a_0\ge a_1-1$. Then $\Delta$ is very admissible.
\end{prop}

\begin{ex}\label{ex:very admissible polygons}
    For any $d,l,k>0$ and $n\ge 0$, the triangle with vertices $(0,0),(d,0),(0,d)$ and the trapezoid with vertices $(0,0),(0,k),(l,k),(l+kn,0)$ clearly satisfy the assumptions of the proposition and, therefore, are very admissible. These are the polygons of the projective plane $\PP^2$ and the Hirzebruch surfaces $\FF_n$ with all possible polarizations. Since a toric blow-up of a surface corresponds to cutting the corners of the corresponding polygon, the polygons of toric blow-ups of $\PP^2$ and $\FF_n$ at the zero-dimensional orbits (including all toric del Pezzo surfaces) equipped with arbitrary polarizations also satisfy the assumptions of the proposition, with the unique exception of $\PP^1\times\PP^1$ blown up at all four zero-dimensional orbits. The polygons of the latter surface can also be shown to be very admissible using a similar but more involved gluing argument than in the proof of Proposition~\ref{prop:example of admissible}. However, the property of being very admissible does not follow from the proposition itself in this case.

    \begin{figure}[ht]   
        \begin{tikzpicture}[x=0.8pt,y=0.8pt,yscale=-1,xscale=1]
            \import{./}{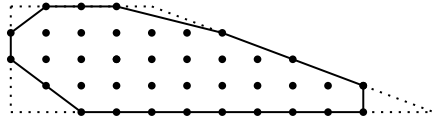}
        \end{tikzpicture}   
    \caption{A polygon of the Hirzebruch surface $\FF_2$ blown up at the four zero-dimensional orbits}\label{fig:polygons}
    \end{figure}
\end{ex}

\begin{rem}
    The polygons of surfaces obtained by further toric blow-ups of the examples above may fail to be very admissible and even $h$-transverse. Indeed, any smooth toric surface can be obtained by repeated toric blow-ups of $\PP^2$ or $\FF_n$, and there are smooth toric surfaces that are not $h$-transverse for any choice of coordinates \cite[Example 3.3]{CHT22}.
\end{rem}

\begin{figure}[ht]   
    \begin{tikzpicture}[x=0.8pt,y=0.8pt,yscale=-.8,xscale=.8]
        \import{./}{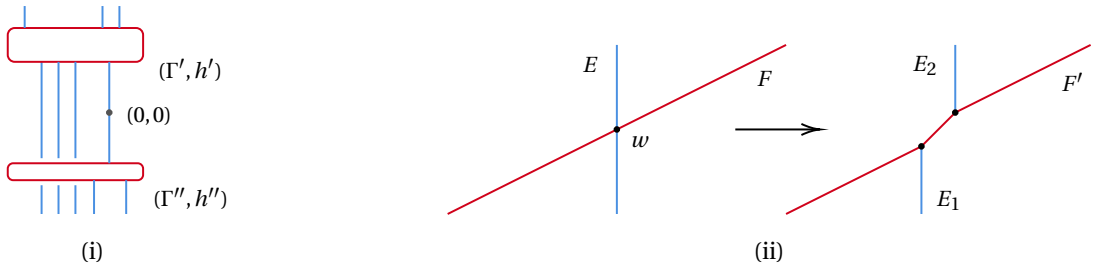}
    \end{tikzpicture}   
    \caption{An illustration to the proof of Proposition~\ref{prop:example of admissible}: (i) The gluing of the parametrized tropical curves $\Gamma'$ and $\Gamma''$, (ii) The perturbation of the curve $(\Gamma_0,h)$.}\label{fig:gluing}
\end{figure}

\begin{ex}
    By \cite[Proposition~4.1.1]{Nod15}, there exist two families of $h$-transverse polarized weighted projective planes: $\PP(1,1,d)$ and $\PP(k,l,k+l)$. The polygons of the first family are the triangles with vertices $(0,0),(0,k),(kd,0)$ for $k,d>0$. These triangles clearly satisfy the assumptions of the proposition, and therefore are very admissible. The $h$-transverse polygons of the second family have no horizontal side (unless $k=l=1$), and therefore are not very admissible.
\end{ex}
We now prove Proposition~\ref{prop:example of admissible}. Recall that a lattice triangle is primitive if it contains no integral points except the vertices, and a \textit{primitive triangulation} of $\Delta$ is a subdivision consisting of primitive triangles.
\begin{proof}[Proof of Proposition~\ref{prop:example of admissible}]
    Let us first prove that $\Delta$ is 0-very admissible. We proceed by induction on the height $\tth$ of $\Delta$ and construct a unimodular ssfd curve of genus 0 with disjoint floors. If $\tth=1$, then the assertion is obvious since any curve dual to a primitive triangulation of $\Delta$ satisfies these properties. For the induction step, consider the intersection $\Delta'$ of $\Delta$ with the half-space given by $y\ge r_0+1$. Since $\Delta$ is convex, the assumptions of the proposition hold true for $\Delta'$, and therefore there exists a unimodular ssfd curve $(\Gamma',h')$ of genus zero with disjoint floors and reduced degree $\nabla'$ dual to $\Delta'$. Furthermore, $(0,-1)\in\nabla'$. Pick a leg $l'$ of $\Gamma'$ whose slope is $(0,-1)$.
    
    Consider the polygon $\Delta''$ of height one, whose non-horizontal sides are parallel to the non-horizontal sides of $\Delta\cap \left(\RR\times [r_0,r_0+1]\right)$, the top side has length one, and the bottom side has length $a_0-a_1+1$, which is non-negative by the assumption of the proposition. Pick any unimodular curve $(\Gamma'',h'')$ of reduced degree $\nabla''$ dual to $\Delta''$, and let $l''$ be the leg of $\Gamma''$ of slope $(0,1)$. After appropriate shifts of $h'$ and $h''$, we may assume that $h''(l'')$ and $h'(l')$ belong to the $y$-axis, $h''$ maps the vertices of $\Gamma''$ sufficiently low in the lower half-plane, and $h'$ maps the vertices of $\Gamma'$ sufficiently high in the upper half-plane. Then the floors of $\Gamma'\cup\Gamma''$ are disjoint. Indeed, by construction, any intersection point of two floors must be an intersection between non-vertical legs of $\Gamma'$ and $\Gamma''$, but such intersections do not exist by the construction and the convexity of $\Delta$.

    To construct the desired curve $(\Gamma,h)$ of degree $\nabla$, it remains to cut the legs $l'$ and $l''$ at the preimages of the origin, and glue the curves $\Gamma'$ and $\Gamma''$ along these cuts. Then $h'$ and $h''$ glue to a parametrization of the curve $\Gamma$, and the resulting curve is a unimodular ssfd curve of genus zero with disjoint floors and degree $\nabla$; see Figure~\ref{fig:gluing}(i).

    To complete the proof, it remains to prove that if $1\le g\le |\Delta^\circ\cap\ZZ^2|$ and there exists a unimodular ssfd curve of genus $g-1$ with disjoint floors and degree $\nabla$, then there exists also such a curve of genus $g$. Let $(\Gamma,h)$ be such a curve of genus $g-1$. We claim that it has at least one self-intersection point. Indeed, if there are no self-intersection points, then $h$ is injective and the image $h(\Gamma)\subset \RR^2$ is a tri-valent graph as $(\Gamma,h)$ is simple. It follows that $(\Gamma,h)$ is dual to a primitive triangulation of $\Delta$. In particular, all inner integral points of $\Delta$ are vertices of the triangulation; and therefore, the curve $(\Gamma,h)$ has genus $|\Delta^\circ\cap\ZZ^2|>g-1$, which is a contradiction.  Pick a self-intersection point $p$. Since the floors of $\Gamma$ are disjoint, $p$ is the intersection of a floor $F$ with an elevator $E$. Consider the curve $(\Gamma_0, h)$ obtained from $\Gamma$ by gluing the preimages of $p$ into a 4-valent vertex $w$, and its perturbation that splits $w$ into a pair of 3-valent vertices such that the new edge is not contracted, see Figure~\ref{fig:gluing}(ii) for an illustration. This perturbation is the desired unimodular, genus-$g$ ssfd curve with disjoint floors and degree $\nabla$. 
\end{proof}

\begin{figure}[ht]   
    \begin{tikzpicture}[x=0.8pt,y=0.8pt,yscale=-0.7,xscale=0.7]
        \import{./}{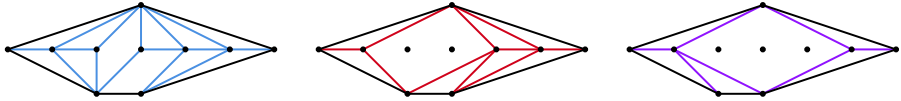}
    \end{tikzpicture}   
    \caption{The polygon in the picture is $g$-admissible if and only if $g\ne 1,2,3$; and is $g$-very admissible if and only if $g=4,5$. In genus 3, it admits a unimodular but not stretched floor decomposed curve. In genera 1 and 2 it admits no unimodular curves, but does admit non-stretched multiplicity-one curves. The dual subdivisions in genera 3, 2, and 1 are depicted in blue, red, and purple respectively.}
    \label{fig:admissible}
\end{figure}

\begin{rem}
    The last part of the argument fails for admissible polygons. In fact, there exist 0-admissible polygons that are not $g$-admissible for some $0<g\le |\Delta^\circ\cap\ZZ^2|$ as the example on Figure~\ref{fig:admissible} shows. 
\end{rem}

\subsection{Zariski's theorem and adjacency of Severi varieties of different genera}\label{subsec:zarANDadj}

Next, we address the geometry of a general curve parametrized by the Severi variety $V_{g, \Delta}^{\irr}$. The following is a generalization of Zariski's theorem about nodality of a general plane curve of given degree and genus. 
\begin{thm}\label{thm:zariski} 
    Suppose $\Delta$ is a $g$-very admissible polygon. Then, over any algebraically closed field $k$, for a general $[C]\in V_{g, \Delta}^{\irr}$, the curve $C$ is at-worst-nodal.
\end{thm}

\begin{cor}\label{cor:adjacencies}
    Suppose $\Delta$ is $(g-1)$-very admissible, and $0<g\le|\Delta^\circ\cap\ZZ^2|$. Then $V^\irr_{g-1,\Delta}\subset \overline V^\irr_{g,\Delta}.$ In particular, if $\Delta$ is very admissible, then 
    $\overline V^\irr_{0,\Delta}\subset \overline V^\irr_{1,\Delta}\subset\dots\subset \overline V^\irr_{g,\Delta}.$
\end{cor}
    
\begin{proof}
    Since $g$-very admissibility is stronger than $g$-admissibility, the Severi variety $V_{g-1, \Delta}^{\irr}$ is irreducible by Theorem~\ref{thm:stongB}. Furthermore, for a general $[C]\in V^\irr_{g-1,\Delta}$, the curve $C$ is at-worst-nodal by Theorem~\ref{thm:zariski}. Since nodes can be smoothed independently (see, e.g., \cite[Lemma~4.5]{Tyo13}), we can construct a partial smoothing of $C$ in which a single node is smoothed. It follows that $V^\irr_{g-1,\Delta}\subset \overline V^\irr_{g,\Delta}$, as claimed.
\end{proof}

\begin{proof}[Proof of Theorem~\ref{thm:zariski}]
    By Theorem~\ref{thm:stongB}, the Severi variety $V_{g, \Delta}^{\irr}$ is irreducible. Therefore, it is sufficient to construct a single nodal curve $[C]\in V_{g, \Delta}^{\irr}$. Furthermore, since being at-worst-nodal is an open condition \cite[Lemma~0DSC]{stacks-project}, it is enough to do so over any algebraically closed field extension of $k$. Thus, we may assume that the ground field is the algebraic closure of the field of Laurent power series $F:=k(\!(\varpi)\!)$. The theorem now follows from the following lemma.
\end{proof}

\begin{lem}\label{lem:nodal curve and unimodular}
    Let $\Delta\subset M_\RR$ be a lattice polygon. Let $[C]\in V^\irr_{g,\Delta}$ be a point such that $C$ intersects the toric boundary of $S_\Delta$ transversely, and the parametrized tropical curve $[\trop(C)]\in M_{g,\nabla}^\trop$ is unimodular. Then $C$ is an at-worst-nodal curve.
\end{lem}
The lemma is an easy consequence of the well-known description of the degeneration of the pair $(C,S_\Delta)$ controlled by the tropicaliztion $\trop(C)$; see, e.g., \cite{Shu05,NS06,Tyo12}. If $\trop(C)$ is unimodular, the degenerated embedded curve is nodal, which implies the nodality of $C$.
\begin{proof}
     Without loss of generality, we may assume that the curve $C$ and its stable model are defined over the discretely valued field $F$.
    
    Let $\sum_{m\in \Delta\cap M}c_mx^m=0$ be an equation defining the curve $C$. Consider the convex hull $\hat{\Delta}$ of the set $\{(m,\nu(c_m))\,|\,m\in \Delta\cap M\; \mathrm{and}\; c_m\ne 0\}\subset\ZZ^3\subset\RR^3$. Since $C$ contains no zero-dimensional orbits, $c_m\ne 0$ for any vertex $m$ of $\Delta$. Thus, $\hat{\Delta}$ induces a piecewise linear function $\xi\: \Delta\to\RR$ given by 
    $$\xi\: m\mapsto \min\{z\,|\,(m,z)\in \hat{\Delta}\}.$$ 
    The linearity domains $\Delta_1,\dotsc,\Delta_r$ of this function provide a subdivision $\subdivision$ of $\Delta$. By Legendre duality, $\subdivision$ is the dual subdivision of $h(\Gamma)$, cf. \cite[\S~2.2]{Shu05}. Now, by the unimodularity of $\trop(C)$, any $\Delta_i$ is either a primitive triangle or a primitive parallelogram.

    Consider the polyhedron $\widetilde{\Delta}:=\{(m,z)\,|\,z\ge \xi(m)\}$. It defines a toric $\ZZ$-scheme $X$ equipped with a line $\CL$ and a toric projective morphism $X\to \AAA^1_\ZZ$. Let $\Spec(F^0)\to \AAA^1_\ZZ=\Spec(\ZZ[t])$ be the morphism given by $t\mapsto\varpi$, and set
    \[\left(S^0,\CL^0\right):=(X,\CL)\times_{\AAA^1_\ZZ}\Spec(F^0).\]
    Then the general fiber of $\left(S^0,\CL^0\right)\to \Spec(F^0)$ is canonically isomorphic to the polarized $F$-scheme $\left(S_\Delta,\CL_\Delta\right)$, and the central fiber is the reducible polarized $\widetilde{F}$-surface $\bigcup_i \left(S_{\Delta_i},\CL_{\Delta_i}\right)$, where for any $i,j$, the surfaces $\left(S_{\Delta_i},\CL_{\Delta_i}\right)$ and $\left(S_{\Delta_j},\CL_{\Delta_j}\right)$ are glued along the invariant closed subset corresponding to $\Delta_i\cap\Delta_j$, cf. \cite[\S~2.3]{Shu05}. Notice that this is the same integral model of $\left(S_\Delta,\CL_\Delta\right)$ as constructed in \cite[\S~3]{Tyo12} in dual terms, i.e., in terms of the tropicalization $\trop(C)$. See also \cite{GS15} for a general theory of toric varieties over rank-one valuation rings.
    
    By construction, the zero locus $\overline{C}$ of the section $\sum_{m\in \Delta\cap M}c_mx^m\in H^0\left(S^0,\CL^0\right)$ is the closure of $C\subset S_\Delta$. Therefore, the intersection $\overline{C}\cap S_{\Delta_i}$ is the zero locus of the section $\sum_{m\in \Delta_i\cap M}c_mx^m$ of $\CL_{\Delta_i}$ for any $i$. In particular, it contains no zero-dimensional orbits of $S_{\Delta_i}$. The morphism $\overline{C}\to\Spec(F^0)$ is flat since $\overline{C}$ is a relative hypersurface in $S^0\to\Spec(F^0)$. Therefore, by \cite[Lemma~0DSC]{stacks-project}, it remains to show that the central fiber of $\overline{C}\to\Spec(F^0)$ is a nodal curve. But the latter is clear, since for any $i$, the polygon $\Delta_i$ is a primitive triangle or a primitive parallelogram, and therefore $\left(S_{\Delta_i},\CL_{\Delta_i}\right)\cong \left(\PP^2,\CO_{\PP^2}(1)\right)$ or $\left(S_{\Delta_i},\CL_{\Delta_i}\right)\cong \left(\PP^1\times\PP^1,\CO_{\PP^1\times\PP^1}(1,1)\right)$. In both cases, the curve $\overline{C}\cap S_{\Delta_i}$ is at worst nodal and intersects the union of other components of the central fiber transversely, which completes the proof.
\end{proof}

\section{Irreducibility of Hurwitz spaces}\label{sec:hurwitz}

In this section, we use the irreducibility of Severi varieties of $\PP^1 \times \PP^1$ to deduce the irreducibility of Hurwitz spaces $\CH_{g,d}$ of all covers (Theorem~\ref{thm:irreducibility of degree d cover}) and Hurwitz spaces $H_{g,d}$ of simple covers (Theorem A). We begin with formal definitions of simple covers, and a lemma showing that the space of simple covers is open in the space of all covers.

\subsection{Coverings} \label{subsec:coverings}

By a $d$-sheeted covering of $\PP^1$ over an algebraically closed field $K$, we mean a finite morphism $\varphi \colon C \to \PP^1$ of degree $d$, where $C$ is a smooth proper curve over $K$. 

\begin{rem}
    Our two main references for coverings are \cite{Ful69} and \cite{LVW24}. In the former, (families of) coverings are defined as finite flat maps; in the latter, as finite locally free maps. The two notions coincide if the base is locally Noetherian (cf. \cite[Tag 02K9]{stacks-project}). In any case, a finite map $\varphi \colon C \to \PP^1$ over an algebraically closed field $K$ as above is automatically flat by \cite[Theorem~23.1]{Mat89}, and therefore locally free by \cite[Tag 02KB]{stacks-project}. 

\end{rem}

In \cite{Ful69}, Fulton introduced the notion of simple $d$-sheeted coverings. If the characteristic of the base field is not 2, simple coverings can be described as follows \cite[Theorem 5.6]{Ful69}.

\begin{defn} 
    Assume that $\mathrm{char}(K) \neq 2$. A {\em simple} $d$-sheeted covering of genus $g$ is a $d$-sheeted covering $\varphi \colon C \to \PP^1$ such that $\varphi$ is separable, $C$ is an irreducible curve of genus $g$, and over each branch point of $\varphi$ there is a single point with non-trivial ramification of index 2.
\end{defn}

In characteristic 2, simple $d$-sheeted coverings in the sense of Fulton do not exist. Either the local degree at a ramification point $p \in C$ is greater than 2, or the covering is not tame. In both cases, it is not simple by \cite[Theorem 2.3]{Ful69}. We modify the definition in this case as follows.

\begin{defn} \label{def:simple char 2}
    Assume that $\mathrm{char}(K)=2$. By a {\em simple} $d$-sheeted covering of genus $g$ we mean a $d$-sheeted covering $\varphi \colon C \to \PP^1$ such that $\varphi$ is separable, $C$ is an irreducible curve of genus $g$, and over each branch point of $\varphi$ there is a single ramification point at which $\varphi$ has local degree 2 and the sheaf of relative differentials has length 2. 
\end{defn}

By definition, the local degree of a simple covering at a ramification point is always 2. If the characteristic is different from 2, then such coverings are tame, and all ramification indices are 2. In characteristic 2, simple coverings are not tame. Since the total length of the sheaf of relative differentials of a separable $d$-sheeted covering of genus $g$ is $2g - 2 + 2d$, it follows that simple coverings have precisely $2g - 2 + 2d$ ramification points if the characteristic is different from 2, and $g - 1 + d$ ramification points if the characteristic is 2. The simplicity of a covering is an open condition, as the following lemma shows.

\begin{lem}\label{lem:openness of simply ramified}
    Let $B$ be a $K$-scheme and $\cC$ a smooth proper relative curve over $B$ with connected geometric fibers. Let $\varphi\colon\cC\rightarrow \PP^1_B$ be a finite locally free morphism of degree $d$. Then the locus $B'$ of points $b\in B$ for which $\cC_b\rightarrow \PP^1$ is simple is open.
\end{lem}
\begin{proof}
    First, the locus $B''\subset B$ of points $b$ such that $\cC_b\rightarrow \PP^1$ is separable is open. Indeed, consider the exact sequence
    $$\varphi^*\Omega_{\PP^1_B/B}\rightarrow \Omega_{\cC/B}\rightarrow \Omega_{\cC/\PP^1_B}\rightarrow 0.$$
    The support $\supp(\Omega_{\cC/\PP^1_B})$ of $\Omega_{\cC/\PP^1_B}$ in $\cC$ is the locus where $\varphi^*\Omega_{\PP^1_B/B}\rightarrow \Omega_{\cC/B}$ is not surjective, hence closed in $\cC$, and proper over $B$. Then $B''$ is the locus of $b$ such that $\supp\left(\Omega_{\cC/\PP^1_B}\right)_b$ has dimension 0. Thus, $B''$ is open in $B$ according to the semi-continuity of fiber dimensions. Set $n:=2g-2+2d$, and consider the branch divisor map $\rho\colon B''\rightarrow{\rm Sym}^n(\PP^1)=\PP^n$ sending $b\in B''$ to the branch divisor of $\cC_b\rightarrow \PP^1$. 
    
    If $\mathrm{char}(K)\neq 2$, then $B'$ is the preimage under $\rho$ of the complement of the union of diagonals in ${\rm Sym}^n(\PP^1)$. Thus, $B'$ is open in $B''$, and hence open in $B$. If $\mathrm{char}(K)=2$, then the image of $B''$ is contained in the sublocus $S$ of even divisors in ${\rm Sym}^n(\PP^1)$. Let $B^*$ be the preimage of the open subset of $S$ consisting of even divisors with support of size $n/2$. Then $B^*$ is open in $B''$, and hence also in $B$. Finally, the locus of points $p\in\PP^1_{B^*}$ such that $\varphi^{-1}(p)$ has at most $d-2$ connected components is closed \cite[Lemma 0BUI]{stacks-project}. Thus, its image $R$ under $\PP^1 _{B^*} \rightarrow B^*$ is closed too, and therefore $B'=B^*\setminus R$ is open.
\end{proof}

\subsection{The Hurwitz spaces of all coverings} 

Recently, Landesman, Vakil and Wood introduced the (algebraic) stack $\CH'_{g,d}$ of $d$-sheeted covers of $\PP^1$ \cite[Definition 5.1]{LVW24}. That is, 
$\CH'_{g,d}$ is the category fibered in groups over $K$-schemes, whose $B$-points over a $K$-scheme $B$ consist of a smooth and proper relative curve $\cC$ over $B$ of arithmetic genus $g$ and a finite locally free map $\cC\rightarrow \PP^1_B$ of degree $d$. According to \cite[12.2.4(vi)]{grothendieck1966elements}, there is an open substack $\CH_{g,d}$ of $\CH'_{g,d}$, whose $B$-points $\cC\rightarrow \PP^1_B$ have geometrically irreducible fibers.

\begin{thm}\label{thm:irreducibility of degree d cover}
    The stack $\CH_{g,d}$ over an algebraically closed field is non-empty and irreducible for any $g\ge 0$ and $d>1$.
\end{thm}
\begin{proof}
    To see that $\CH_{g,d}$ is irreducible, it suffices to find an irreducible variety dominating $\CH_{g,d}$. Pick $c\gg 1$ coprime to $d$ such that every curve $C$ of genus $g$ admits a degree-$c$ covering $C \to \PP^1$. Let $\Delta$ be a rectangle with horizontal sides of length $c$ and vertical sides of length $d$. Then the polarized toric surface $(S_\Delta,\CL_\Delta)$ is isomorphic to $\left(\PP^1\times\PP^1,\CO_{\PP^1\times\PP^1}(c,d)\right)$, and by Theorem~B, the Severi variety $V^\irr_{g,\Delta}$ is irreducible.

    Let $\mathcal C \to V^\irr_{g,\Delta}$ be the universal curve over $V^\irr_{g,\Delta}$. Arguing as in the proof of Lemma~\ref{lem:surjectivity into nice cones}, there is an alteration $V \to V^\irr_{g,\Delta}$, with $V$ irreducible, such that the pullback $\CC_V$ of the universal curve is equinormalizable over $V$, i.e., the fibers of the normalization of the total space $\mathcal C^\nu_V \to V$ are all smooth of genus $g$. Furthermore, the composition $\mathcal C^\nu_V \to \PP^1 \times \PP^1\xrightarrow{\pi}\mathbb P^1$ is a family of $d$-sheeted coverings of $\PP^1$ over $V$, where $\pi$ is the projection to the first factor. This family induces a map $\varphi \colon  V \rightarrow \CH_{g,d}.$ Thus, $\CH_{g,d}$ is non-empty, and it remains to prove that $\varphi$ is dominant. 
    
    Given a point $[(C, f_1 \colon C \to \PP^1)]\in \CH_{g,d}$, pick a covering $f_2 \colon C \to \PP^1$ of degree $c$. Then the map $(f_1,f_2)\colon C \to \PP^1 \times \PP^1$ is birational onto its image. Indeed, if not, then there is a factorization $C \to C' \to \PP^1 \times \PP^1$, where the map $C \to C'$ has degree $k > 1$. But then $f_1$ and $f_2$ must factor through $C \to C'$, and therefore $k>1$ is a common divisor of $c$ and $d$, which is a contradiction since $c$ and $d$ are coprime. Thus, $[(C,f_1)]\in \varphi(V)$, and we are done.
\end{proof}

\begin{rem}
    According to \cite[12.2.4(vi)]{grothendieck1966elements}, we can decompose $\CH_{g,d}'$ into a disjoint union of open substacks, where each substack specifies the number of geometric connected components of each fiber $\cC_b$ of $\cC\rightarrow B$, the genus $g'$ of each component of $\cC_b$, and the degree $d'$ of the restriction of the covering of $\PP^1$ to the component. Since all such substacks are dominated by products of $\CH_{g',d'}$'s, it follows from Theorem~\ref{thm:irreducibility of degree d cover} that all such substacks are irreducible. 
\end{rem}

\subsection{The Hurwitz spaces of simple coverings}

In this subsection, we finally prove Theorem~A. 

\begin{ThmA}
    The Hurwitz space $H_{g,d}$ of simple $d$-sheeted coverings of $\PP^1$ over an algebraically closed field is non-empty and irreducible for all $g\ge 0$ and $d>1$.
\end{ThmA}

\begin{proof}
    The irreducibility of $H_{g,d}$ follows from Theorem~\ref{thm:irreducibility of degree d cover} and the fact that $H_{g,d}\subseteq |\CH_{g,d}|$ is open by Lemma~\ref{lem:openness of simply ramified}. It remains to show that there exists a simple $d$-sheeted covering for any degree $d>1$ and genus $g\ge 0$. As in the proof of Theorem~\ref{thm:irreducibility of degree d cover}, let $\Delta$ be the rectangle with horizontal sides of length $c$ and vertical sides of length $d$, where $c\gg 1$. Then $(S_\Delta,\CL_\Delta)$ is isomorphic to $\left(\PP^1\times\PP^1,\CO_{\PP^1\times\PP^1}(c,d)\right)$. Denote by $\pi\colon\PP^1\times\PP^1\to\PP^1$ the projection to the first factor.

    Let $L_0$ be a general curve of bidegree $(c,1)$, and let $L_1,\dotsc, L_{d-1}$ be general curves of bidegree $(0,1)$. Set $C_0:=\bigcup_{i=0}^{d-1} L_i$. Then $C_0$ is a nodal curve of bidegree $(c,d)$ with $c(d-1)$ nodes. Furthermore, $\pi$ is injective on the set of nodes. Pick a collection $\CP\subset C_0$ of $d+g-1$ nodes such that $\CP$ contains at least one node of $L_0\cap L_i$ for all $i\ge 1$. Such a collection exists because $c\gg 1$. Note that the $c(d-1)$ nodes of $C_0$ can be smoothed independently; see, e.g., \cite[Lemma~4.5]{Tyo13}. The desired simple covering is then constructed by considering a general smoothing of the nodes of $\CP$ in a partial resolution of $C_0$.

    Consider a one-parameter deformation $C_t$ of $C_0$ over $K[\![t]\!]$ that smooths the nodes of $\CP$ and preserves the other nodes. The normalization $C'$ of its geometric generic fiber is a smooth irreducible curve of genus $g$, and the projection $\pi'\: C'\to \PP^1$ is a separable covering of degree $d$. Because $\pi'$ is separable, the Riemann-Hurwitz formula applies. Thus, the sum of the lengths of the sheaf of relative differentials at the ramification points is given by 
    \begin{equation}\label{eq: totalram}
        2g(C')-2+2d=2|\CP|.
    \end{equation}
    
    We claim that $\pi'\: C'\to \PP^1$ is a simple covering. To prove this in $\mathrm{char}(K)\neq 2$ (resp., $\mathrm{char}(K)=2$), it suffices to show that for each node $p\in\CP$, there exist two ramification points (resp., one ramification point) of local degree 2 specializing to $p$, and if $\mathrm{char}(K)\neq 2$, the corresponding branch points are distinct. Indeed, because the sum of the lengths of the sheaf of relative differentials at those ramification points is at least $2|\CP|$, equation \eqref{eq: totalram} implies that $\pi'$ has no additional ramification, and the length of the sheaf of relative differentials is 1 (resp., 2) at each ramification point. Furthermore, because $\pi$ is injective on $\CP$, all branch points are distinct. Therefore, the covering $\pi'\: C'\to \PP^1$ is simple. 
    
    Let $p\in\CP$ be any point. Then the germ of $C_t$ at $p$ is given by $y(y-x) = \lambda$ for some non-zero $\lambda\in tK[\![t]\!]$, and $\pi$ is the projection onto the $x$-axis. Therefore, $\pi'$ has two (resp., one) ramification points of local degree 2 specializing to $p$, given locally by $(x,y)=\pm (2\sqrt{-\lambda}, \sqrt{-\lambda})$. If $\mathrm{char}(K) \neq 2$, the corresponding branch points are indeed distinct, which completes the proof. 
\end{proof}

\bibliographystyle{amsalpha}
\bibliography{1}
\end{document}